\documentclass[12pt]{article}
\usepackage[margin=1in]{geometry}
\usepackage[T1]{fontenc}
\usepackage{lmodern,microtype}
\usepackage{amsmath,amssymb,amsthm,mathtools,bm}
\usepackage{enumitem,booktabs,tabularx,array}
\usepackage[round,authoryear]{natbib}
\usepackage{xcolor}
\usepackage{pgfplots}
\pgfplotsset{compat=1.18}
\newenvironment{keywords}{\par\noindent\textbf{Keywords:} }{\par}
\setlist{leftmargin=1.5em,itemsep=0.2em,topsep=0.3em}
\theoremstyle{plain}
\newtheorem{theorem}{Theorem}
\newtheorem{corollary}{Corollary}
\newtheorem{proposition}{Proposition}
\newtheorem{lemma}{Lemma}
\newtheorem{definition}{Definition}
\theoremstyle{definition}
\newtheorem{assumption}{Assumption}
\newtheorem{example}{Example}
\theoremstyle{remark}
\newtheorem{remark}{Remark}

\newcommand{\HH}{\mathcal H}\newcommand{\CC}{\mathcal C}\newcommand{\KK}{\mathcal K}\newcommand{\QQ}{\mathcal Q}
\newcommand{\RR}{\mathbb R}\newcommand{\ZZ}{\mathbb Z}\newcommand{\one}{\mathbf 1}
\newcommand{\ip}[2]{\langle #1,#2\rangle}\newcommand{\norm}[1]{\lVert #1\rVert}
\newcommand{\cconv}{\overline{\operatorname{conv}}}
\DeclareMathOperator{\lin}{lin}\DeclareMathOperator{\ran}{ran}\DeclareMathOperator{\supp}{supp}
\DeclareMathOperator{\spn}{span}\DeclareMathOperator{\Id}{Id}\DeclareMathOperator{\ncap}{cap}
\DeclareMathOperator*{\essinf}{ess\,inf}\DeclareMathOperator*{\esssup}{ess\,sup}
\newcommand{\Expect}[1]{\mathbb{E}(#1)}\newcommand{\ExpectBracket}[1]{\mathbb{E}[#1]}\newcommand{\ExpectSquared}[1]{\mathbb{E}\{(#1)^2\}}
\newcommand{\ProbMeasure}{\mathbb{P}}\newcommand{\ProbOf}[1]{\mathbb{P}(#1)}
\newcommand{\VarOf}[1]{\operatorname{var}(#1)}\newcommand{\TraceOf}[1]{\operatorname{tr}(#1)}
\newcommand{\Differential}[1]{\,\mathrm{d}#1}
\newcommand{\AlmostSure}{almost surely}\newcommand{\AlmostSureEnd}{almost surely.}\newcommand{\AlmostSureAdjective}{almost sure}
\newcommand{\AlmostEverywhere}{almost everywhere}\newcommand{\AlmostEverywhereEnd}{almost everywhere.}\newcommand{\AlmostEvery}{almost every}
\newcommand{\IndependentIdentical}{independent and identically distributed}\newcommand{\FirstOrderAR}{first-order autoregressive}
\newcommand{\TextAbbr}[2]{#2}
\newcommand{\tbl}[2]{\caption{#1}#2}\newcommand{\ReferenceTitle}[1]{}\newcommand{\AppendixPart}[1]{\subsubsection*{#1}}
\usepackage[hidelinks]{hyperref}
\title{Admissibility of Linear Time Series Models under Cone Constraints}
\author{Won-Ki Seo\\[0.5ex] School of Economics, University of Sydney}
\date{}
\begin{document}
\maketitle
\begin{abstract}
Time series analysis extends beyond vector-valued observations to evolving probability densities, quantile
functions, and covariance operators. These objects carry intrinsic constraints, raising questions about the
admissibility of additive linear models, which have long served as a foundation for time series modelling and
theory. We give a necessary and sufficient support condition for a stationary Hilbert-space linear process to
remain in a closed convex cone almost surely. Innovations with unbounded lower tails, including Gaussian
innovations, confine randomness to the cone's largest linear subspace. For fractional filters with convergent
normalized cumulative coefficients, the leading stochastic directions obey the same restriction, even with
bounded innovations. For real-line quantiles, variation under unbounded lower tails and leading fractional accumulation can occur
only through location shifts.  For cone-valued processes, partial sums of centred stationary increments remain
bounded in mean square when evaluated by any linear functional that is nonnegative on the cone, in contrast to
the growing variance of a random walk. We derive an exact slope criterion for locally absolutely continuous
quantiles and construct nondegenerate additive models with bounded innovations for quantiles, densities, and
covariance matrices. These findings clarify the scope of ordinary additive modelling under structural
constraints.
\end{abstract}


\section{Introduction}\label{sec_intro}
Linear models form a common foundation for conventional and functional time series analysis. Describing
observations through the effects of current and past innovations yields tractable models of temporal dependence
and a basis for prediction, estimation, and  inference \citep{BrockwellDavis1991,PhillipsSolo1992}.
Hilbert-space extensions carry this approach to functional autoregression and linear processes
\citep{Bosq2000}, representations of stochastic trends and cointegration
\citep{BeareSeoSeo2017,BeareSeo2020,FranchiParuolo2020}, and inference on the dimension of common trends
\citep{NielsenSeoSeong2023,NielsenSeoSeongForthcoming}. This combination of modelling simplicity and theoretical
depth makes linear models a natural benchmark for new kinds of time series data.

One such extension is to time series of distributions, where the objective is to follow changes in dispersion,
asymmetry, or tail behaviour beyond a single scalar summary. Densities and quantile functions represent these
observations as curves, making the functional linear framework a natural starting point. Covariance operators
provide another example of objects embedded in a linear space. Yet these representations must obey constraints:
densities are nonnegative and normalized, quantiles are nondecreasing, and covariance operators are positive
semidefinite. Applying the linear framework therefore requires determining whether its realizations remain
valid objects and how much variation an admissible model can accommodate.

Quantile functions are used as functional observations in distributional regression
\citep{YangEtAl2020,GhosalEtAl2023} and, when recorded over successive periods, in autoregressive models
\citep{ZhangKokoszkaPetersen2022,ZhuMuller2023}. For univariate distributions with finite second moments, the
quantile map is an isometric embedding of the $2$-Wasserstein space into $L^2(0,1)$ \citep{PanaretosZemel2020}.
Its image is the closed convex cone of functions with a nondecreasing version. The Hilbert-space representation
thus retains a monotonicity constraint, which remains after centring at a reference quantile. Moreover, this
cone has empty interior in $L^2(0,1)$, so an arbitrarily small perturbation in the Hilbert norm can disrupt
quantile ordering. A linear representation must therefore be accompanied by conditions that preserve
monotonicity.

Preserving these constraints also depends on how shock effects persist. Shocks that preserve the constraints
individually can violate them when combined across lags. In a quantile autoregression developed below, greater
persistence reduces the maximum admissible shock magnitude. Fractional filters used to model long memory raise
a related issue because their cumulative coefficients can grow without bound. The support and persistence of
shocks are thus joint modelling choices. A useful characterization must explain both the restrictions and the
scope for nondegenerate dynamics.

Under suitable, method-specific assumptions, existing approaches accommodate distributional constraints through
transformations to Hilbert-space coordinates \citep{PetersenMuller2016}, modified linear operations in Bayes
Hilbert spaces \citep{SeoBeare2019}, or Wasserstein tangent coordinates and optimal transport maps
\citep{ZhangKokoszkaPetersen2022,ZhuMuller2023}. In particular, \citet{ZhuMuller2023} formulate autoregressive
dynamics using compositions of transport maps and a form of scalar multiplication that preserves monotonicity.
Our results provide a complementary perspective by characterizing when ordinary additive models are admissible
and how innovation support and persistence restrict the variation they can represent. Suitable support and
filter conditions also permit nondegenerate additive models with bounded innovations. Structural constraints are inherent to distribution-valued data and must be respected regardless of the modelling framework.

We study a real separable Hilbert space $\HH$ and a closed convex cone $\CC\subseteq\HH$. Let $\mu\in\HH$ be a
fixed mean, let $\{\varepsilon_t\}$ be centred independent and identically distributed $\HH$-valued
innovations, and let $a_j:\HH\to\HH$ be bounded linear coefficients. We consider the stationary linear process
\begin{equation}\label{eq_linear}
X_t=\mu+\sum_{j\ge0}a_j\varepsilon_{t-j},\qquad t\in\ZZ,
\end{equation}
whose series converges in mean square. We characterize when this model remains in $\CC$ and which forms and
magnitudes of stochastic variation it can represent.  Appendix~\ref{app_proofs} proves the main results. The Supplementary Material gives further results and examples for each section, with their proofs.
\section{Admissibility through innovation support}\label{sec_support}
\subsection{Geometry and the linear process}
\TextAbbr{We write a.s.\ for almost surely and a.e.\ for almost everywhere. Unless another measure is specified, a.e.\ refers to Lebesgue measure.}{Unless another measure is specified, ``almost everywhere'' refers to Lebesgue measure.}

Write $\ip{\cdot}{\cdot}$ and $\norm{\cdot}$ for the inner product and norm on $\HH$, and $\mathcal B(\HH)$ for
the space of bounded linear operators on $\HH$, with operator norm $\norm{\cdot}_{\mathrm{op}}$. For
$D\subseteq\HH$, $\spn D$, $D^\perp$, and $\overline D$ denote its linear span, orthogonal complement, and norm
closure, respectively. For an operator $A$, $A^*$, $\ran A$, and $\ker A$ denote its adjoint, range, and
kernel. The identity operator is denoted by $\Id$, and $(x\otimes y)z=\ip{y}{z}x$ for $x,y,z\in\HH$. The
principal examples use $L^2(0,1)$, the Hilbert space of Lebesgue square-integrable real functions on $(0,1)$,
with functions equal \AlmostEverywhere{} identified and inner product
$\ip{f}{g}=\int_0^1f(u)g(u)\Differential{u}$. We abbreviate it to $L^2$ when no confusion arises. For a
measurable real function $f$ on $(0,1)$, write $\norm{f}_\infty=\esssup_{u\in(0,1)}|f(u)|$, allowing the value
$\infty$. If $f$ is continuous on $[0,1]$, this equals $\sup_{u\in[0,1]}|f(u)|$.

The dual cone is $\CC^*=\{h\in\HH:\ip{h}{x}\ge0\text{ for all }x\in\CC\}$. The bipolar theorem for closed
convex cones \citep[Section~6.3]{BauschkeCombettes2017} gives $x\in\CC$ if and only if $\ip{h}{x}\ge0$ for
every $h\in\CC^*$. Applying this characterization to both $x$ and $-x$ yields
\begin{equation}\label{eq_lineality}
\lin(\CC)=(\CC^*)^\perp=(\overline{\spn}\,\CC^*)^\perp.
\end{equation}
A cone with $\lin(\CC)=\{0\}$ is called pointed. We define its \emph{noise capacity} as
$\ncap(\CC)=\dim\lin(\CC)$, allowing the value $\infty$. Its role depends on the innovation or accumulation
conditions established below; it does not bound the dimension of bounded admissible randomness.

\begin{assumption}[Linear process]\label{ass_linear}
The innovations $\{\varepsilon_t\}_{t\in\ZZ}$ are independent and identically distributed $\HH$-valued random
elements on a probability space $(\Omega,\mathcal F,\ProbMeasure)$, with $\Expect{\varepsilon_0}=0$ and
$\Expect{\norm{\varepsilon_0}^2}<\infty$. Their covariance operator is
$\Sigma=\ExpectBracket{\varepsilon_0\otimes\varepsilon_0}$, and $\Sigma^{1/2}$ denotes its nonnegative square
root. The coefficients $a_j\in\mathcal B(\HH)$ satisfy $\sum_{j\ge0}\Expect{\norm{a_j\varepsilon_0}^2}<\infty$,
and $\mu\in\HH$.
\end{assumption}
 For each fixed $t$, the series in \eqref{eq_linear} converges in mean square
and \AlmostSureEnd{}  The process is strictly stationary and has
mean $\mu$. This coefficient condition is weaker than absolute summability in operator norm and includes
square-summable scalar filters; see \citet{Bosq2000} for the standard linear process framework.

Let $\cconv$ denote the norm closure of the convex hull. Without a subscript, $\essinf$ and $\esssup$ are taken under the probability measure. Writing $\supp\varepsilon_0$ for the topological support of the innovation law, define
$\KK=\cconv(\supp\varepsilon_0)$ and $\kappa(g)=\inf_{x\in\KK}\ip{g}{x}=\essinf\ip{g}{\varepsilon_0}$ for $g\in\HH$.
The map $\kappa$, which we call the lower support functional, assigns to each $g$ the essential lower bound of
the scalar innovation projection $\ip{g}{\varepsilon_0}$. Its expression as an infimum over $\KK$ follows from
continuity and the definition of support. The innovation mean belongs to $\KK$, so $0\in\KK$ and
$\kappa(g)\le0$. 
\subsection{Support criterion for admissibility}
\begin{theorem}[Admissibility characterization]\label{thm_support}
Under Assumption~\ref{ass_linear}, the following statements are equivalent:
\begin{enumerate}[label=(\roman*),ref=(\roman*)]
\item $X_t\in\CC$ \AlmostSure{} for every $t\in\ZZ$.
\item\label{item_support_budget} For every $h\in\CC^*$, $\ip{h}{\mu}+\sum_{j\ge0}\kappa(a_j^*h)\ge0$.
\item\label{item_finite_support} For every integer $m\ge0$ and every $x_0,\ldots,x_m\in\KK$, $\mu+\sum_{j=0}^m
a_jx_j\in\CC$.
\end{enumerate}
Moreover, $\essinf\ip{h}{X_t}=\ip{h}{\mu}+\sum_{j\ge0}\kappa(a_j^*h)$ for every $h\in\HH$. Every series of
lower support values $\kappa(a_j^*h)$ is interpreted in $[-\infty,0]$; any $-\infty$ term or divergence to
$-\infty$ violates condition~\ref{item_support_budget}.
\end{theorem}
Condition~\ref{item_support_budget} requires the mean projection $\ip{h}{\mu}$ to cover the total worst
downward contribution $-\sum_j\kappa(a_j^*h)$ of the lagged innovation terms, for every $h\in\CC^*$. Here
$\kappa(a_j^*h)=\essinf\ip{h}{a_j\varepsilon_0}$ is the essential lower bound of the lag-$j$ contribution in
that direction. Condition~\ref{item_finite_support} provides a convenient test for admissibility in the
constructions below: it involves only deterministic finite sums, even though $X_t$ is defined by an infinite
series. Checking these sums suffices for admissibility because, under Assumption~\ref{ass_linear}, the realized
partial sums converge \AlmostSure{} to $X_t$ and $\CC$ is closed. Convergence for arbitrary infinite sequences
from $\KK$ is not required.

For $\HH=\RR$ and $\CC=[0,\infty)$, condition~\ref{item_support_budget} reads $\mu+\sum_{j\ge0}\essinf(a_j\varepsilon_0)\ge0$, which holds, for example, for first-order autoregressions with a nonnegative coefficient and positive square-integrable innovations; \citet{DavisMcCormick1989} study estimation in autoregressions with positive innovations.

Section~\ref{sec_support_criterion_s} of the Supplementary Material derives two consequences of Theorem~\ref{thm_support}: admissibility requires $\mu\in\CC$, and if $\ip{h}{\mu}=0$ for some $h\in\CC^*$, every lagged term $\ip{h}{a_j\varepsilon_{t-j}}$ vanishes \AlmostSure{}, so a mean in $\lin(\CC)$ forces $X_t-\mu\in\lin(\CC)$. The same section discusses the role of independence: by Remark~\ref{rem_independence}, with dependent innovations the support criterion is sufficient but not necessary.

\section{Restrictions from support and persistence}\label{sec_collapse}
We examine how innovation support and the persistence of shock effects restrict admissible linear dynamics.
\subsection{Unbounded lower tails}
\begin{definition}\label{def_unbounded}
A centred innovation with covariance $\Sigma$ and lower support functional $\kappa$ has \emph{unbounded lower
tails} if $\Sigma^{1/2}g\ne0$ implies $\kappa(g)=-\infty$ for every $g\in\HH$.
\end{definition}
Definition~\ref{def_unbounded} requires each scalar projection $\ip{g}{\varepsilon_0}$ with positive variance $\norm{\Sigma^{1/2}g}^2$ to be unbounded below: $\ProbOf{\ip{g}{\varepsilon_0}<-M}>0$ for every $M>0$. Projections with zero variance vanish \AlmostSure{}, so degenerate covariances are allowed. Section~\ref{sec_unbounded_s} of the Supplementary Material shows this and two further points: every projection with positive variance is then also unbounded above, and the condition is implied neither by an unbounded innovation norm nor by a Gaussian component that is degenerate in a direction with positive total innovation variance. 

A centred Gaussian random element satisfies the condition: each
scalar projection $\ip{g}{\varepsilon_0}$ is either identically zero or has support $\RR$ \citep{Bogachev1998}.
The following theorem shows that,
under this tail condition, admissible stochastic variation is confined to the cone's lineality space.

\begin{theorem}[Collapse to lineality]\label{thm_collapse}
Suppose Assumption~\ref{ass_linear} holds and the innovation has unbounded lower tails. Then $X_t\in\CC$
\AlmostSure{} if and only if $\mu\in\CC$ and
\begin{equation}\label{eq_range}
\ran(a_j\Sigma^{1/2})\subseteq\lin(\CC)\qquad\text{for every }j\ge0.
\end{equation}
In that case $a_j\varepsilon_{t-j}\in\lin(\CC)$ \AlmostSure{}, $X_t\in\mu+\lin(\CC)$ \AlmostSure{}, and the
covariance operator of $X_t$ has range contained in $\lin(\CC)$. If the cone is pointed, $X_t=\mu$
\AlmostSureEnd{}
\end{theorem}
Nonzero rescaling of the innovations cannot restore admissibility when \eqref{eq_range} fails. A scalar
projection $\ip{g}{\varepsilon_0}$ with positive variance remains unbounded below after rescaling, however
small its variance becomes. Conversely, any centred linear process taking values in lineality is admissible
after translation by a mean in the cone. Thus $\dim\lin(\CC)$ is the exact available dimension for this
innovation class, and not merely an upper bound without an attainable counterpart.
Example~\ref{ex_one_sided} in Section~\ref{sec_unbounded_s} of the Supplementary Material shows that unbounded support alone does not suffice: centred exponential innovations, unbounded above but bounded below, give a nondegenerate process in $[0,\infty)$.

\subsection{Admissibility of persistent linear processes}
The linear process in Assumption~\ref{ass_linear} can also exhibit persistence, including long-range
dependence, since absolute summability of the coefficients $a_j$ is not required. We now examine how
admissibility restricts persistence without an unbounded-tail assumption. The cumulative coefficient
$F_m=\sum_{j=0}^m a_j$ records the response to a single innovation accumulated through lag $m$. We allow the
coefficient operators to vary across lags and impose a limit only on the normalized cumulative response.
\begin{theorem}[Restriction on cumulative responses]\label{thm_accumulation}
Suppose Assumption~\ref{ass_linear} holds and $X_t\in\CC$ \AlmostSureEnd{} If positive scalars $r_m\to\infty$
and an operator $\widetilde F_\infty\in\mathcal B(\HH)$ satisfy $r_m^{-1}F_m\to\widetilde F_\infty$ in operator
norm, then
\begin{equation}\label{eq_accumulation_range}
\ran(\widetilde F_\infty\Sigma^{1/2})\subseteq\lin(\CC).
\end{equation}
In particular, suppose that, for some $0<d<1/2$ and $B\in\mathcal B(\HH)$, the coefficients satisfy
\begin{equation}\label{eq_fractional_filter}
a_j=j^{d-1}B+R_j\quad(j\ge1),\qquad
\frac{\norm{R_j}_{\mathrm{op}}}{j^{d-1}}\longrightarrow0
\quad\text{as }j\to\infty.
\end{equation}
Then $m^{-d}F_m\to B/d$ in operator norm and $\ran(B\Sigma^{1/2})\subseteq\lin(\CC)$. These conclusions hold
for every innovation law satisfying Assumption~\ref{ass_linear}, including bounded innovations.
\end{theorem}

In the fractional case \eqref{eq_fractional_filter}, the restriction concerns $\ran(B\Sigma^{1/2})$: for
quantile functions of real-valued random variables, this range contains only constant functions; on a pointed
cone, it is zero. Absolutely summable components may still produce admissible variation outside $\lin(\CC)$.
The location shift can therefore follow fractional dynamics while shape variation remains bounded, as
illustrated after Example~\ref{ex_ar}.

Remark~\ref{rem_scalar} in Section~\ref{sec_scalar_s} of the Supplementary Material, with details in Section~\ref{sec_scalar_details}, shows that, for scalar multiples $a_j=b_jB$ with $\sum_j|b_j|=\infty$, as in the usual scalar fractional filter and in functional local Whittle methods \citep{LiRobinsonShang2021,SeoShang2024}, admissibility confines the whole centred process to $\lin(\CC)$, even when $\sum_{j\le m}b_j$ stays bounded. The reason is that, for each $h\in\CC^*$ with $\VarOf{\ip{h}{B\varepsilon_0}}>0$, the lag-$j$ lower support value $\kappa(b_jB^*h)$ is at most $-c_h|b_j|$ for some $c_h>0$; these values sum to $-\infty$, so condition~\ref{item_support_budget} of Theorem~\ref{thm_support} fails.

\section{Processes with stationary increments}\label{sec_increments}
Consider $Q_t=Q_0+\sum_{s=1}^tu_s\in\CC$ with centred stationary increments $u_s$. Recall that $h\in\CC^*$
implies $\ip{h}{Q_t}\ge0$. The following theorem shows that, for every such $h$, the scalar partial sum
$\sum_{s=1}^t\ip{h}{u_s}=\ip{h}{Q_t-Q_0}$ is the difference of two terms of a square-integrable stationary
sequence. Its second moment is therefore bounded uniformly in $t$.
\begin{theorem}[Cone-valued levels with stationary increments]\label{thm_increments}
Let $\{u_t\}_{t\ge1}$ be a strictly stationary $\HH$-valued sequence with $\Expect{u_1}=0$ and
$\Expect{\norm{u_1}^2}<\infty$. Let $Q_0\in\CC$ \AlmostSure{} with $\Expect{\norm{Q_0}^2}<\infty$, and suppose
that $Q_t=Q_0+\sum_{s=1}^t u_s\in\CC$ \AlmostSure{} for every $t\ge1$.
For every $h\in\CC^*$, there is a strictly stationary square-integrable scalar sequence $\{M_t^h\}_{t\ge0}$,
with each $M_t^h$ measurable with respect to $\{u_{t+s}\}_{s\ge1}$, such that
\begin{equation}\label{eq_projected_difference}
\sum_{s=1}^t\ip{h}{u_s}=M_0^h-M_t^h,\qquad
\ExpectSquared{M_0^h}\le\Expect{\ip{h}{Q_0}^2}.
\end{equation}
Consequently,
\begin{equation}\label{eq_uniform_projection}
\sup_{t\ge1}\Expect{\ip{h}{Q_t-Q_0}^2}
\le4\Expect{\ip{h}{Q_0}^2}<\infty.
\end{equation}
\end{theorem}
For each fixed $h\in\CC^*$, the sequence in \eqref{eq_projected_difference} can be chosen as
$M_t^h=\inf_{n\ge0}\sum_{s=1}^n\ip{h}{u_{t+s}}=\inf_{s\ge t}\ip{h}{Q_s}-\ip{h}{Q_t}$,
with the empty sum equal to zero. Thus $M_t^h$ is the infimum of $\ip{h}{Q_s-Q_t}$ over $s\ge t$. Since
$h\in\CC^*$ and $Q_s\in\CC$, we have $\ip{h}{Q_s}\ge0$ and hence $-\ip{h}{Q_t}\le M_t^h\le0$. The theorem
assumes no independence between $Q_0$ and $\{u_t\}_{t\ge1}$, and none among the increments.
For stationary sequences, such representations exist when the partial sums have uniformly bounded second moments \citep[Lemma~5 and the proof of Theorem~2]{Browder1958}; here the cone constraint gives the representation directly.
\begin{corollary}[Limits of accumulated increments]\label{cor_limits}
Under the conditions of Theorem~\ref{thm_increments}, put $S_T=Q_T-Q_0=\sum_{s=1}^T u_s$.
\begin{enumerate}[label=(\roman*)]
\item If $\{r_T\}$ is a positive deterministic sequence with $r_T\to\infty$ and $S_T/r_T$ converges weakly in
$\HH$ to $G$, then $G\in\lin(\CC)$ \AlmostSureEnd{}
\item If a bounded nonnegative self-adjoint operator $\Lambda$ satisfies
\begin{equation}\label{eq_lrv}
\ip{h}{\Lambda h}=\lim_{T\to\infty}T^{-1}\Expect{\ip{h}{S_T}^2}
\quad\text{for every }h\in\HH,
\end{equation}
then $\ran\Lambda\subseteq\lin(\CC)$ and $\operatorname{rank}\Lambda\le\ncap(\CC)$${}=\dim\lin(\CC)$.
\end{enumerate}
\end{corollary}
Section~\ref{sec_increments_conditions} of the Supplementary Material shows that the centring condition $\Expect{u_1}=0$ cannot be dropped, and that Theorem~\ref{thm_increments} restricts accumulated increments, not the dependence of the levels: for any strictly stationary $\CC$-valued process with $\Expect{\norm{Q_0}^2}<\infty$, $Q_T-Q_0$ is bounded in mean square, however slowly the serial covariance of $Q_t$ decays. Section~\ref{sec_increments_corollary} explains how Corollary~\ref{cor_limits} follows from Theorem~\ref{thm_increments} and shows that the limits in Corollary~\ref{cor_limits}(ii) exist, for example, under absolutely summable lag covariances. It also explains that the scalar bounds give $\sup_T\Expect{\norm{PS_T}^2}<\infty$, for the orthogonal projection $P$ onto $\lin(\CC)^\perp$, when $\lin(\CC)^\perp$ is finite dimensional, but do not by themselves give it otherwise. Section~\ref{sec_density_normalization} treats random densities whose mean may vary over time. By Remark~\ref{rem_density}, unit mass alone forces normalized limits of changes in centred random densities to vanish, giving the restriction of \citet{Beare2017} for any square-integrable densities, without stationary increments or linearity.

\section{Geometry of the principal examples}\label{sec_geometry}
\noindent Quantile functions of univariate distributions with finite second moment form the closed
convex cone \begin{equation}\label{eqquantilecone}\QQ=\{q\in L^2(0,1):q\text{ has a nondecreasing version}\}.\end{equation} The quantile map is an isometry from univariate $2$-Wasserstein space into $L^2(0,1)$ \citep{PanaretosZemel2020}.
Write $\one(u)=1$ for the constant function on $(0,1)$. A function and its negative
are both nondecreasing precisely when the function is constant, so $\lin(\QQ)=\spn\{\one\}$. Quantiles of
distributions supported on $[0,\infty)$ are nonnegative, which makes their cone pointed. Centred quantiles
satisfy $\int_0^1q(u)\Differential{u}=0$, which also removes the constant direction. Both cases therefore have
zero lineality, as do the nonnegative cone of $L^2(0,1)$, which contains square-integrable densities, and the positive semidefinite cone of self-adjoint Hilbert--Schmidt operators. In all these cones, Theorem~\ref{thm_collapse} makes additive processes with unbounded lower tails deterministic; for densities, nonnegativity alone forces this.

Section~\ref{sec_geometry_s} of the Supplementary Material gives the details for densities and positive semidefinite operators, and treats convex functions, whose lineality space consists of affine functions, so under the same condition only the intercept and slope vary. Independent Gaussian random walks in the intercept and slope show that the bound $\operatorname{rank}\Lambda\le\ncap(\CC)=2$ in Corollary~\ref{cor_limits}(ii) is then attained. Monotonicity leaves only the intercept, and adding nonnegativity, as for Lorenz curves, gives the same degeneracy as for densities. Table~\ref{tab_geometry} summarizes these cases. By Proposition~\ref{prop_affine}, an affine restriction $BX_t=b$, such as unit mass or a fixed mean, holds \AlmostSure{} if and only if $B\mu=b$ and $Ba_j\Sigma^{1/2}=0$ for every $j$. Under unbounded lower tails it confines $X_t-\mu$ to $\lin(\CC)\cap\ker B$, and in Theorem~\ref{thm_accumulation} it gives $\ran(\widetilde F_\infty\Sigma^{1/2})\subseteq\lin(\CC)\cap\ker B$.

\section{Admissible linear dynamics of quantile functions}\label{sec_quantile}
Admissible linear models for quantile functions are time series models for distributions in $2$-Wasserstein geometry \citep{ZhangKokoszkaPetersen2022}. When the innovations have unbounded lower tails, Theorem~\ref{thm_collapse} confines quantile variation to location shifts; Gaussian innovations are one example. Theorem~\ref{thm_accumulation} does the same for the leading fractional term in \eqref{eq_fractional_filter}, for any innovation law. These results might suggest that admissible quantile models can vary only within $\lin(\QQ)=\spn\{\one\}$, that is, only in location. Monotonicity alone does not force this. Bounded innovations permit changes in shape when
their magnitude and the coefficients $a_j$ satisfy suitable bounds, even though $\QQ$ has empty interior in $L^2(0,1)$. We first derive an exact slope condition
for admissibility under local absolute continuity, without assuming bounded innovations. We then use it to construct such models with one, several, or infinitely many random shape
directions.

Throughout this section, $\HH=L^2(0,1)$ and $\CC=\QQ$ (see \eqref{eqquantilecone}). Recall that $\KK=\cconv(\supp\varepsilon_0)$ is the
closed convex innovation support from Section~\ref{sec_support}, and $\mu=\Expect{X_t}$ is the mean quantile
function.
\subsection{An exact slope condition}
To characterize admissible linear dynamics through derivatives, we assume that $\mu$ and each function $a_jx$,
$x\in\KK$, are locally absolutely continuous. This means absolute continuity on every closed subinterval of
$(0,1)$. For such a function $h$, $h'$ denotes its derivative, defined \AlmostEverywhereEnd{}
\begin{assumption}[Local absolute continuity]\label{ass_slope}
The mean $\mu\in\QQ$ has a locally absolutely continuous version on $(0,1)$. For every $j\ge0$ and every
$x\in\KK$, $a_jx$ has a locally absolutely continuous version on $(0,1)$.
\end{assumption}
For each lag define, up to Lebesgue null sets, the measurable function $s_j:(0,1)\to[0,\infty]$ by
$s_j=\mathop{\mathrm{ess\,sup}}_{x\in\KK}\{-(a_jx)'\}$.
Since $0\in\KK$, $s_j\ge0$.

\begin{theorem}[Necessary and sufficient slope condition]\label{thm_slope}
Under Assumptions~\ref{ass_linear} and~\ref{ass_slope},
\begin{equation}\label{eq_exact_slope}
X_t\in\QQ\text{ \AlmostSure{}}
\quad\Longleftrightarrow\quad
\mu'(u)\ge\sum_{j\ge0}s_j(u)\quad\text{for \AlmostEvery{} }u\in(0,1).
\end{equation}
\end{theorem}
In \eqref{eq_exact_slope}, $s_j(u)$ measures how much the term $a_jx$ can reduce the slope as $x$ ranges over
$\KK$. The sum over lags must be no greater than the mean slope $\mu'(u)$. The innovation support determines
the possible slopes at each lag, and the coefficients determine how their contributions add across lags.  Theorem~\ref{thm_slope} complements the sufficient derivative bounds
of \citet{ZhangKokoszkaPetersen2022} by providing an exact admissibility criterion for general operator filters
under the stated regularity.

Section~\ref{sec_slope_interp} of the Supplementary Material interprets \eqref{eq_exact_slope}: if $\mu$ is the quantile function of a positive density $f_\mu$, then $\mu'(u)=1/f_\mu\{\mu(u)\}$ wherever the inverse-function relation applies, so larger slope perturbations are allowed where $f_\mu$ is small; on a flat interval of $\mu$, which represents a point mass, $X_t$ must stay flat but can shift in location. It also notes that the slope criterion needs local absolute continuity, since otherwise a nonnegative derivative \AlmostEverywhere{} need not imply monotonicity; Theorem~\ref{thm_support} does not need it.

\subsection{Constructing admissible linear processes}
\begin{figure}[t]
\centering
\begin{tikzpicture}
\begin{axis}[width=0.36\linewidth,height=3.8cm,xmin=0,xmax=1,ymin=-0.4,ymax=1.4,domain=0:1,samples=150,xlabel={$u$},ylabel={$q_\pm(u)$},title={Admissible: $\rho=0.6,\ \lambda=0.3$},title style={font=\footnotesize},tick label style={font=\footnotesize},label style={font=\footnotesize},axis lines=left,xtick={0,0.5,1},ytick={0,0.5,1}]
\addplot[gray,densely dashed] {x};
\addplot[blue!70!black,thick] {x+0.75*x*(1-x)};
\addplot[blue!70!black,thick] {x-0.75*x*(1-x)};
\end{axis}
\begin{axis}[at={(0.43\linewidth,0)},anchor=south west,width=0.36\linewidth,height=3.8cm,xmin=0,xmax=1,ymin=-0.4,ymax=1.4,domain=0:1,samples=150,xlabel={$u$},title={Inadmissible: $\rho=0.6,\ \lambda=1.2$},title style={font=\footnotesize},tick label style={font=\footnotesize},label style={font=\footnotesize},axis lines=left,xtick={0,0.5,1},ytick={0,0.5,1}]
\addplot[gray,densely dashed] {x};
\addplot[red!70!black,thick] {x+3*x*(1-x)};
\addplot[red!70!black,thick] {x-3*x*(1-x)};
\end{axis}
\end{tikzpicture}
\caption{Support extremes $q_\pm$ in Example~\ref{ex_ar} with $\rho=0.6$; the dashed line is the mean $\mu(u)=u$. Left, $\lambda=0.3$: both are nondecreasing. Right, $\lambda=1.2$: $q_-$ decreases on $(0,1/3)$ and $q_+$ on $(2/3,1)$. These are deterministic bounds, not simulated paths.}\label{fig_boundary}
\end{figure}
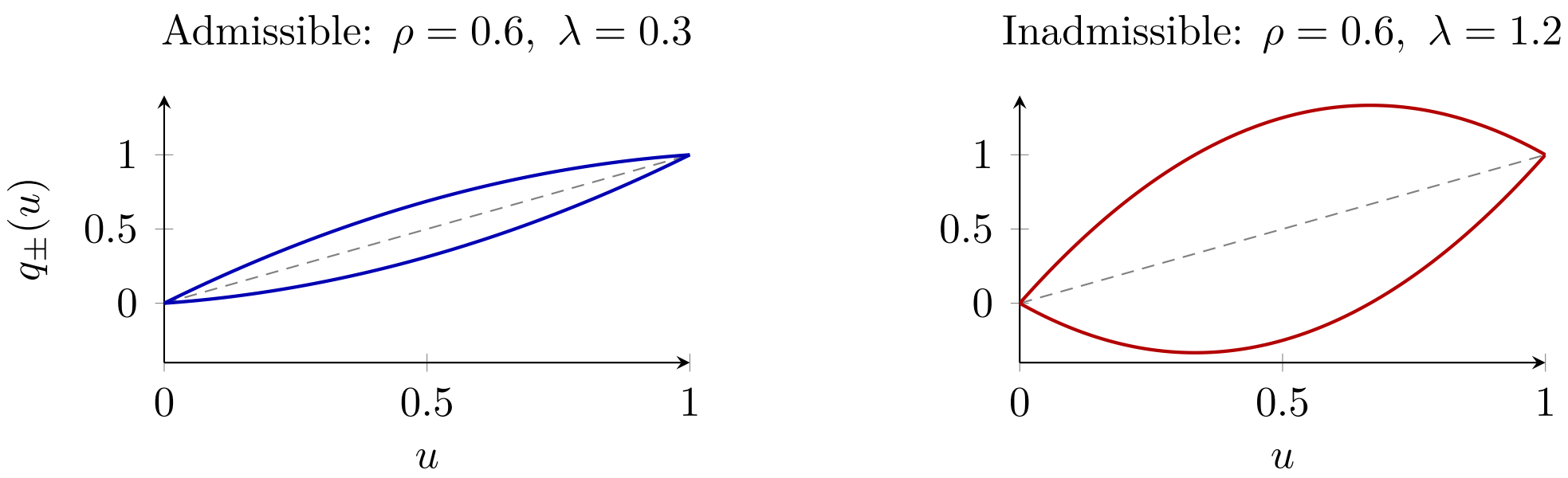

We construct admissible linear quantile processes with bounded innovations. Example~\ref{ex_ar} gives a model
with one shape direction, for which Theorem~\ref{thm_slope} yields an exact admissibility condition.
Proposition~\ref{prop_infinite} in Section~\ref{sec_many_directions} of the Supplementary Material extends the construction to finitely or infinitely many shape directions under
a sufficient bound on their combined downward slopes.
\begin{example}[One shape direction]\label{ex_ar}
Let $\mu(u)=u$, $v(u)=u(1-u)$, and let $\eta_t$ be independent and identically distributed and centred, with positive variance and $|\eta_t|\le1$. For $|\rho|<1$ and
$\lambda>0$, consider
$X_t(u)=u+Z_tu(1-u)$ with $Z_t=\lambda\sum_{j\ge0}\rho^j\eta_{t-j}$.

A realized function is nondecreasing precisely if $|Z_t|\le1$. Since $|Z_t|\le\lambda/(1-|\rho|)$, the process is quantile-valued \AlmostSure{} whenever $\lambda\le1-|\rho|$. This uses only the bound on $\eta_t$, not its distribution. If, in addition, $\essinf\eta_0=-1$ and $\esssup\eta_0=1$, as for Rademacher variables or uniform variables on $[-1,1]$, then $s_j(u)=\lambda|\rho|^j|1-2u|$, and Theorem~\ref{thm_slope} shows that $\lambda\le1-|\rho|$ is also necessary. The essential extrema of $Z_t$ are then $\pm\lambda/(1-|\rho|)$. Its support need not fill the interval between these extrema. If
$\lambda>1-|\rho|$, a positive-probability support event produces a nonmonotone function.
Figure~\ref{fig_boundary} plots the support extremes $q_\pm(u)=u\pm\{\lambda/(1-|\rho|)\}u(1-u)$ of $X_t$ on either side of this boundary.
\end{example}
Every realization lies between $q_-$ and $q_+$ and, under the endpoint condition, comes arbitrarily close to each with positive probability. In the right panel of Figure~\ref{fig_boundary}, monotonicity is lost near the end points, where the slope $|v'(u)|=|1-2u|$ of the shape direction is largest relative to the mean slope $\mu'(u)=1$. The boundary $\lambda\le1-|\rho|$ also quantifies how persistence restricts shape variation. 
Under the endpoint condition in Example~\ref{ex_ar}, $\VarOf{Z_t}=\lambda^2\VarOf{\eta_0}/(1-\rho^2)\le\lambda^2/(1-\rho^2)$, so admissibility with $0\le\rho<1$ implies $\VarOf{Z_t}\le(1-\rho)/(1+\rho)$. Thus persistence approaching one forces the shape variance to zero. Section~\ref{sec_ar_example} of the Supplementary Material discusses the example further: how the sign of $Z_t$ moves probability mass in the represented distributions, a direct argument that the extremes of $Z_t$ are $\pm\lambda/(1-|\rho|)$ under the endpoint condition, and why the support of $Z_t$, a Cantor set for Rademacher $\eta_t$ and $|\rho|<1/2$ \citep{PeresSchlagSolomyak2000}, does not affect admissibility.

Any centred square-integrable scalar process $L_t$ can be added as a location shift: $X_t+L_t\one$ remains a
quantile under $\lambda\le1-|\rho|$. In particular, $L_t$ may be a Gaussian process. Examples~3.2--3.3 of
\citet{ZhangKokoszkaPetersen2022} already allow unbounded location innovations alongside bounded shape
variation under a common scalar autoregressive filter. Here, location and shape can follow different filters.
For example, take $L_t=\sum_{j\ge0}(j+1)^{d-1}\zeta_{t-j}$, where $0<d<1/2$ and $\{\zeta_t\}$ are independent and identically distributed, centred and square integrable, also independent of $\{\eta_t\}$. The location shift then follows fractional dynamics,
while the shape coefficient $Z_t$ retains its \FirstOrderAR{} dynamics. This illustrates how
Theorem~\ref{thm_accumulation} permits fractional accumulation in location together with bounded shape
variation. 

Proposition~\ref{prop_infinite} in Section~\ref{sec_many_directions} of the Supplementary Material extends Example~\ref{ex_ar} to countably many shape directions $v_k$ with coefficients bounded by one, with $a_j=\lambda\rho^j\Id$, under $\lambda\sum_{k}\norm{v_k'}_\infty\le c(1-|\rho|)$ when $\mu'\ge c>0$. With $\mu(u)=u$, $\varepsilon_t=\sum_{k}\eta_{t,k}v_k$ and $v_k(u)=2^{-k}(2\pi k)^{-1}\sin(2\pi ku)$, we have $\sum_k\norm{v_k'}_\infty=1$, and this gives admissibility for $\lambda\le1-|\rho|$, with innovation and process covariance operators of infinite rank. The same section shows that, with the standard normal quantile as the mean, $\mu'\ge\sqrt{2\pi}$ and $\lambda\sum_k\norm{v_k'}_\infty\le\sqrt{2\pi}(1-|\rho|)$ suffices.

These examples separate the size of shape perturbations from their dimension. The mean slope and lagged shock
effects determine how much shape change is admissible; bounded innovations can meet this condition in
infinitely many directions. Location changes remain unrestricted by monotonicity.

\section{Admissible dynamics for general cones}\label{sec_bounded}
As for quantile functions in Section~\ref{sec_quantile}, bounded innovations allow nondegenerate additive models in general cones, even pointed ones. Section~\ref{sec_bounded_s} of the Supplementary Material uses the gauge $p_\mu(v)=\sup_{h\in\CC^*}|\ip{h}{v}|/\ip{h}{\mu}$, with $0/0=0$, which satisfies $p_\mu(v)\le1$ exactly when $\mu\pm v\in\CC$. By Proposition~\ref{prop_gauge}, for $\mu\in\CC$, $\sum_{j\ge0}\sup_{x\in\KK}p_\mu(a_jx)\le1$ suffices for admissibility; this condition bounds only the perturbations $a_jx$, so $\mu$ need not be an interior point of $\CC$. Because it treats each lag separately, it can be conservative, whereas condition~\ref{item_support_budget} of Theorem~\ref{thm_support} is exact. By Examples~\ref{ex_density}--\ref{ex_quantile_general}, derived in Section~\ref{sec_relative_details}, $p_\mu(v)$ equals $\norm{v/\mu}_\infty$ for densities, $\norm{\mu^{-1/2}v\mu^{-1/2}}_{\mathrm{op}}$ for covariance matrices with positive definite mean, and $\norm{v'/\mu'}_\infty$ for quantiles when $\mu$ and $v$ are locally absolutely continuous. For instance, with $Z_t$ as in Example~\ref{ex_ar}, $X_t(u)=1+Z_t(2u-1)$ is a nondegenerate density-valued process whenever $\lambda\le1-|\rho|$, so the tail condition in Theorem~\ref{thm_collapse} is needed. Proposition~\ref{prop_gauge} also covers dependent innovations with values in $\KK$ when the filtered series converges \AlmostSure{}.

\section{Discussion and concluding remarks}\label{sec_discussion}
Admissibility jointly restricts the constraint set, the innovation distribution, and the coefficient operators. Unbounded lower tails confine variation to $\lin(\CC)$, and accumulation does the same for leading fractional terms and for limits of accumulated stationary increments, even with bounded innovations. For quantiles, only location shifts remain; for densities and covariance operators, no randomness remains. Bounded innovations nevertheless permit nondegenerate models, with exact conditions given by Theorems~\ref{thm_support} and~\ref{thm_slope}.

These restrictions matter for model specification. If a functional autoregression for quantile curves is simulated, or used for prediction, with Gaussian innovations that are not almost surely constant functions, Theorem~\ref{thm_collapse} implies that the generated curves are nonmonotone with positive probability, however small their nonconstant part. Monotone rearrangement \citep{ChernozhukovEtAl2010} restores monotonicity, but the rearranged curves no longer follow the additive dynamics. Theorem~\ref{thm_slope} instead shows how to restrict innovation support and persistence so that the model itself remains admissible.

Further work could derive explicit restrictions for general operator coefficients without a nonzero limit of their normalized cumulative sums, and could treat dependent innovations, for which admissibility depends on which innovation values can occur together at different lags.

\section*{Declaration of generative AI use}
ChatGPT Astra 6.0 (OpenAI) was used to check the mathematical arguments, and Claude Opus 5.5 (Anthropic) was used to edit the language. The author reviewed and verified all content and takes full responsibility for it.

\appendix
\renewcommand{\thelemma}{A.\arabic{lemma}}
\section{Proofs of the main results}\label{app_proofs}
We use the innovation support $\KK=\cconv(\supp\varepsilon_0)$, its lower support functional $\kappa(g)=\inf_{x\in\KK}\ip{g}{x}$, and the innovation covariance $\Sigma$ from Section~\ref{sec_support}, and write $\mathbb{L}^2_{\mathbb{R}}$ for the space of real random variables with finite second moment. The functional $\kappa$ is positively homogeneous for positive scalars and superadditive. Its negative is the usual support function of $\KK$ evaluated at $-g$.
\begin{proof}[Proof of Theorem~\ref{thm_support}]
We first prove the essential-infimum identity stated in the final part of the theorem:
$\essinf\ip{h}{X_t}=\ip{h}{\mu}+\sum_{j\ge0}\kappa(a_j^*h)$. We then establish the equivalence of (i), (ii),
and (iii). Throughout the proof, we use the identity $\essinf(Y+Z)=\essinf Y+\essinf Z$ for independent,
\AlmostSure{} finite real random variables $Y$ and $Z$, with essential infima allowed to equal $-\infty$. This
follows by taking infima in the support identity of \citet[Appendix~A, Proposition~2.1]{SteutelVanHarn2004}.

Fix $t\in\ZZ$ and $h\in\HH$, and put $\xi_j=\ip{a_j^*h}{\varepsilon_{t-j}}$,
$\ell_j=\essinf\xi_j=\kappa(a_j^*h)$, and $S=\sum_{j\ge0}\xi_j=\ip{h}{X_t-\mu}$. Assumption~\ref{ass_linear}
ensures that the series converges \AlmostSure{} and in $\mathbb{L}^2_{\mathbb{R}}$. The variables $\xi_j$ are
independent and centred, so $\ell_j\le0$ for every $j$. Consequently, $\sum_{j\ge0}\ell_j$ is well defined as
the limit of its nonincreasing partial sums, with values in $[-\infty,0]$.

For each $m\ge0$, write $S_m=\sum_{j=0}^m\xi_j$ and $R_m=S-S_m=\sum_{j>m}\xi_j$. The variables $S_m$ and $R_m$
are independent because they depend on disjoint sets of innovations. Moreover, convergence in mean square and
centring give $\Expect{R_m}=0$, so $\essinf R_m\le0$. Repeated application of $\essinf(Y+Z)=\essinf Y+\essinf
Z$ to the independent terms in $S_m$ gives $\essinf S_m=\sum_{j=0}^m\ell_j$. Applying the same equality to
$S=S_m+R_m$ then yields $\essinf S=\sum_{j=0}^m\ell_j+\essinf R_m\le\sum_{j=0}^m\ell_j$. Letting $m\to\infty$
gives $\essinf S\le\sum_{j\ge0}\ell_j$. If this sum is $-\infty$, equality follows immediately. Otherwise, on a
probability-one event, $S_m\to S$ and $\xi_j\ge\ell_j$ for every $j$. Taking limits in
$S_m\ge\sum_{j=0}^m\ell_j$ gives $S\ge\sum_{j\ge0}\ell_j$, and hence the reverse inequality for $\essinf S$.
Adding the constant $\ip{h}{\mu}$ proves the essential-infimum identity for $X_t$.

We now prove (i)$\Rightarrow$(ii)$\Rightarrow$(iii)$\Rightarrow$(i).

\emph{(i) $\Rightarrow$ (ii).}
Suppose $X_t\in\CC$ \AlmostSure{} for every $t\in\ZZ$. Fix $h\in\CC^*$ and any $t\in\ZZ$. By the definition of
$\CC^*$, $\ip{h}{X_t}\ge0$ \AlmostSureEnd{} The formula for $\essinf\ip{h}{X_t}$ therefore gives
$\ip{h}{\mu}+\sum_{j\ge0}\kappa(a_j^*h)=\essinf\ip{h}{X_t}\ge0$, which is (ii).

\emph{(ii) $\Rightarrow$ (iii).}
Fix $m\ge0$ and $x_0,\ldots,x_m\in\KK$. For every $h\in\CC^*$, the definition of $\kappa$ gives
$\ip{a_j^*h}{x_j}\ge\kappa(a_j^*h)$. Since every $\kappa(a_j^*h)\le0$, a finite partial sum of these values is
no smaller than their full sum. Thus (ii) implies
\begin{equation*}
\ip{h}{\mu+\sum_{j=0}^m a_jx_j}
\ge\ip{h}{\mu}+\sum_{j=0}^m\kappa(a_j^*h)
\ge\ip{h}{\mu}+\sum_{j\ge0}\kappa(a_j^*h)\ge0.
\end{equation*}
These deterministic inequalities hold for every $h\in\CC^*$. The dual-cone characterization recalled in
Section~\ref{sec_support} therefore gives $\mu+\sum_{j=0}^m a_jx_j\in\CC$, proving (iii).

\emph{(iii) $\Rightarrow$ (i).}
Fix $t\in\ZZ$. Since $\ProbOf{\varepsilon_0\in\KK}=1$, (iii) gives $\mu+\sum_{j=0}^m
a_j\varepsilon_{t-j}\in\CC$ \AlmostSure{} for every $m\ge0$. These sums converge to $X_t$ in mean square, so a
subsequence converges almost surely to $X_t$. Since every partial sum lies in $\CC$ \AlmostSure{} and $\CC$ is
closed, $X_t\in\CC$ \AlmostSureEnd{}
\end{proof}

\begin{proof}[Proof of Theorem~\ref{thm_collapse}]
Suppose first that $X_t\in\CC$ \AlmostSureEnd{} Since $\CC$ is closed and convex, $\mu=\Expect{X_t}\in\CC$.
Condition~\ref{item_support_budget} of Theorem~\ref{thm_support} gives $\kappa(a_j^*h)>-\infty$ for every
$h\in\CC^*$ and every $j\ge0$. Definition~\ref{def_unbounded} thus implies $\Sigma^{1/2}a_j^*h=0$. By
continuity, the bounded operator $\Sigma^{1/2}a_j^*$ vanishes on $\overline{\spn}\CC^*=\lin(\CC)^\perp$, where
the equality follows from \eqref{eq_lineality}. Taking adjoints yields \eqref{eq_range}.

Conversely, suppose $\mu\in\CC$ and \eqref{eq_range} holds. We first show that
$\varepsilon_0\in\overline{\ran\Sigma^{1/2}}=(\ker\Sigma)^\perp$ \AlmostSureEnd{} Choose an orthonormal basis
$\{g_k\}$ of $\ker\Sigma$. For each $k$, $\Expect{\ip{g_k}{\varepsilon_0}^2}=\ip{g_k}{\Sigma g_k}=0$, so
$\ip{g_k}{\varepsilon_0}=0$ \AlmostSureEnd{} A countable intersection makes these equalities simultaneous for
all $k$, implying that $\varepsilon_0$ is orthogonal to $\ker\Sigma$. The same conclusion holds for every
innovation $\varepsilon_t$ by identical distribution.  For any $x\in\overline{\ran\Sigma^{1/2}}$, choose
$z_n\in\HH$ such that $\Sigma^{1/2}z_n\to x$. Equation~\eqref{eq_range} gives $a_j\Sigma^{1/2}z_n\in\lin(\CC)$
for every $n$. Continuity of $a_j$ and closedness of $\lin(\CC)$ thus imply $a_jx\in\lin(\CC)$.

Applying this to the innovations gives $a_j\varepsilon_{t-j}\in\lin(\CC)$ \AlmostSure{} for every $j$. Every
finite sum of these terms thus belongs to $\lin(\CC)$. Passing to an almost surely convergent subsequence and
using closedness of $\lin(\CC)$, we obtain $X_t-\mu\in\lin(\CC)$ and hence $X_t\in\mu+\lin(\CC)\subseteq\CC$
\AlmostSureEnd{} A countable intersection makes these conclusions simultaneous for all $t\in\ZZ$ and $j\ge0$.

Finally, $\ip{h}{X_t-\mu}=0$ \AlmostSure{} for all $h\in\lin(\CC)^\perp$, so the covariance  of $X_t$ vanishes
on $\lin(\CC)^\perp$. Self-adjointness thus places its range in $\lin(\CC)$. If $\CC$ is pointed,
$\lin(\CC)=\{0\}$ and $X_t=\mu$ \AlmostSureEnd{}
\end{proof}

Condition~\ref{item_finite_support} of Theorem~\ref{thm_support} permits the same innovation-support value
$x\in\KK$ at every lag, giving $\mu+F_mx\in\CC$. Dividing by $r_m$ and passing to the limit places $\widetilde
F_\infty x$ in the cone. Since $\widetilde F_\infty\varepsilon_0$ is centred, each $\ip{h}{\widetilde
F_\infty\varepsilon_0}$ with $h\in\CC^*$ is nonnegative with mean zero and must vanish. This proves the
restriction without a limit theorem for process partial sums.

\begin{proof}[Proof of Theorem~\ref{thm_accumulation}]
First suppose $r_m^{-1}F_m\to\widetilde F_\infty$ in operator norm, with $r_m\to\infty$. Recall
$F_m=\sum_{j=0}^m a_j$. For every $x\in\KK$, condition~\ref{item_finite_support} of Theorem~\ref{thm_support},
applied with $x_0=\cdots=x_m=x$, gives $\mu+F_mx\in\CC$. Since $r_m>0$, the cone contains
$r_m^{-1}\mu+r_m^{-1}F_mx$. Closedness implies $\widetilde F_\infty x\in\CC$.

Hence $\widetilde F_\infty\varepsilon_0\in\CC$ \AlmostSureEnd{} For each $h\in\CC^*$, the scalar
$\ip{h}{\widetilde F_\infty\varepsilon_0}$ is nonnegative and has mean zero, so it vanishes \AlmostSureEnd{}
Its variance is $\norm{\Sigma^{1/2}\widetilde F_\infty^*h}^2=0$. Equation~\eqref{eq_lineality} and taking
adjoints yield \eqref{eq_accumulation_range}.

For the fractional case \eqref{eq_fractional_filter}, take $r_m=m^d$ for $m\ge1$. Since $x^{d-1}$ is decreasing
on $[1,\infty)$, $\int_1^{m+1}x^{d-1}\,dx\le\sum_{j=1}^m j^{d-1}\le1+\int_1^m x^{d-1}\,dx$. Evaluating the
integrals and dividing by $m^d$ gives $m^{-d}\sum_{j=1}^m j^{d-1}\to1/d$.
Given $\delta>0$, choose $J\ge1$ such that $\norm{R_j}_{\mathrm{op}}\le\delta j^{d-1}$ for $j\ge J$. The
contribution from $j<J$, divided by $m^d$, tends to zero, while the remaining terms have norm at most $\delta
m^{-d}\sum_{j=J}^m j^{d-1}$. Thus $\limsup_{m\to\infty}\norm{m^{-d}\sum_{j=1}^m R_j}_{\mathrm{op}}\le\delta/d$.
Letting $\delta\downarrow0$ gives $m^{-d}\sum_{j=1}^m R_j\to0$ in operator norm. Since $m^{-d}a_0\to0$, it
follows that $m^{-d}F_m\to B/d$. Applying \eqref{eq_accumulation_range} with $\widetilde F_\infty=B/d$ yields
$\ran(B\Sigma^{1/2})\subseteq\lin(\CC)$, since $1/d\ne0$.
\end{proof}

\begin{proof}[Proof of Theorem~\ref{thm_increments}]
Fix $h\in\CC^*$, write $x_t=\ip{h}{u_t}$, and define $M_t^h=\inf_{n\ge0}\sum_{s=1}^n x_{t+s}$ for $t\ge0$, with
the empty sum equal to zero. The infimum is over a countable set, so $M_t^h$ is measurable with respect to
$\{u_{t+s}\}_{s\ge1}$. Since $Q_n\in\CC$ \AlmostSure{}, $\sum_{s=1}^n x_s=\ip{h}{Q_n-Q_0}\ge-\ip{h}{Q_0}$ \AlmostSure{}
for every $n\ge0$. Hence $-\ip{h}{Q_0}\le M_0^h\le0$, which gives $M_0^h\in\mathbb{L}^2_{\mathbb{R}}$ and
$\ExpectSquared{M_0^h}\le\Expect{\ip{h}{Q_0}^2}$.

For each $t$, $M_t^h$ is obtained by applying the same measurable map
$(y_s)_{s\ge1}\mapsto\inf_{n\ge0}\sum_{s=1}^n y_s$ to $(x_{t+s})_{s\ge1}$. Strict stationarity of $\{x_t\}$
thus implies strict stationarity of $\{M_t^h\}$. Each $M_t^h$ has the same distribution as $M_0^h$ and thus
belongs to $\mathbb{L}^2_{\mathbb{R}}$.

For $t\ge1$, the definition gives $M_{t-1}^h=\min(0,x_t+M_t^h)\le x_t+M_t^h$. The integrable variable
$x_t+M_t^h-M_{t-1}^h$ is nonnegative and has mean zero, because $\Expect{x_t}=0$ and $M_t^h$ and $M_{t-1}^h$
have the same distribution. It is therefore zero \AlmostSureEnd{} Hence $x_t=M_{t-1}^h-M_t^h$, and summing
gives \eqref{eq_projected_difference}. Finally, stationarity gives
$\ExpectSquared{M_0^h-M_t^h}\le4\ExpectSquared{M_0^h}$, proving \eqref{eq_uniform_projection}.
\end{proof}
\begin{proof}[Proof of Corollary~\ref{cor_limits}]
Recall $S_T=\sum_{s=1}^T u_s$. For every fixed $h\in\CC^*$, \eqref{eq_uniform_projection} gives
$\ip{h}{S_T}/r_T\to0$ in $\mathbb{L}^2_{\mathbb{R}}$ for the normalization $r_T\to\infty$ in part (i). If
$S_T/r_T\Rightarrow G$, the continuous mapping theorem applied to $x\mapsto\ip{h}{x}$ gives $\ip{h}{G}=0$
\AlmostSureEnd{} Choose a countable dense subset $\{h_n\}_{n\ge1}$ of $\CC^*$. Taking a countable intersection
gives a probability-one event on which $\ip{h_n}{G}=0$ for every $n$. On this event, continuity of
$h\mapsto\ip{h}{G}$ extends the equality to every $h\in\CC^*$. Since $\ip{h}{G}=\ip{h}{-G}=0$ for every such
$h$, the dual characterization of a closed convex cone implies that both $G$ and $-G$ belong to $\CC$. Hence
$G\in\CC\cap(-\CC)=\lin(\CC)$ \AlmostSureEnd{}

For part (ii), let $\Lambda$ satisfy \eqref{eq_lrv}. The same bound gives $\ip{h}{\Lambda h}=0$ for every
$h\in\CC^*$. Nonnegativity of $\Lambda$ implies $\Lambda h=0$ for every such $h$. The closed linear subspace
$\ker\Lambda$ therefore contains $\overline{\spn}\,\CC^*=\lin(\CC)^\perp$. Self-adjointness implies
$\ran\Lambda\subseteq\lin(\CC)$, from which $\operatorname{rank}\Lambda\le\ncap(\CC)$ immediately follows.
\end{proof}

The proof of Theorem~\ref{thm_slope} uses the following lemma.
\begin{lemma}[Countable representation]\label{lem_essential_family}
Let $(S,\mathcal A,m)$ be a finite measure space, and let $\{g_a:S\to\mathbb R\}_{a\in I}$ be a nonempty family
of measurable functions. There exists a countable set $J\subseteq I$ such that
\begin{equation*}
g:=\mathop{\mathrm{ess\,sup}}_{a\in I}g_a
=\sup_{a\in J}g_a
\qquad m\text{-\AlmostEverywhere{}}
\end{equation*}
\end{lemma}
\begin{proof}[Proof of Lemma~\ref{lem_essential_family}]
The zero-measure case is immediate. Otherwise, normalize $m$ to a probability measure. This preserves null
sets, and the functions $g_a$ are random variables on the resulting probability space. Let $\Phi=\{\max_{a\in
F}g_a:F\subseteq I,\ F\text{ finite and nonempty}\}$. The maximum of two members of $\Phi$ is again a finite
maximum, so $\Phi$ is directed upward \citep[Theorem~A.37]{FollmerSchied2016}. Each $g_a$ belongs to $\Phi$ by
taking $F=\{a\}$. Conversely, if $\psi\ge g_a$ \AlmostEverywhere{} for every $a$, taking the maximum over any
finite $F$ gives $\psi\ge\max_{a\in F}g_a$ \AlmostEverywhereEnd{} Thus $\Phi$ and the original family have the
same \TextAbbr{\AlmostEverywhere{}}{almost-everywhere} upper bounds. By \citet[Appendix~A.5,
Theorem~A.37]{FollmerSchied2016}, their essential supremum $g$ exists, and there is an increasing sequence
$h_n\in\Phi$ with $h_n\uparrow g$ \AlmostEverywhereEnd{} For each $n$, write $h_n=\max_{a\in F_n}g_a$ for a
finite set $F_n\subseteq I$. Then $J=\bigcup_{n\ge1}F_n$ is countable, and $g=\sup_{n\ge1}h_n=\sup_{a\in J}g_a$
\AlmostEverywhereEnd{}
\end{proof}

\begin{proof}[Proof of Theorem~\ref{thm_slope}]
Recall $s_j=\mathop{\mathrm{ess\,sup}}_{x\in\KK}(-(a_jx)')$ from Section~\ref{sec_quantile}. By
Lemma~\ref{lem_essential_family}, for each $j$ we can choose $x_{j,k}\in\KK$, $k\ge1$, with $x_{j,1}=0$, such
that $s_j(u)=\sup_{k\ge1}\{-(a_jx_{j,k})'(u)\}$ \AlmostEverywhereEnd{} For $n\ge1$, put $g_{j,n}(u)=\max_{1\le
k\le n}\{-(a_jx_{j,k})'(u)\}$. Thus $g_{j,n}$ is the largest of the first $n$ selected negative derivatives at
lag $j$, and $0\le g_{j,n}\uparrow s_j$ \AlmostEverywhereEnd{}

Suppose first that $X_t\in\QQ$ \AlmostSureEnd{} Fix $m\ge0$ and $n\ge1$, and consider all choices
$x_j\in\{x_{j,1},\ldots,x_{j,n}\}$ for $j=0,\ldots,m$. Condition~\ref{item_finite_support} of
Theorem~\ref{thm_support} gives $\mu+\sum_{j=0}^m a_jx_j\in\QQ$ for every such choice. Its locally absolutely
continuous version is nondecreasing, so $\mu'\ge\sum_{j=0}^m(-(a_jx_j)')$ \AlmostEverywhereEnd{}  For fixed
$m,n$, there are at most $n^{m+1}$ choices of $(x_0,\ldots,x_m)$. Thus only countably many derivative
inequalities are involved as $m,n$ vary. After discarding a Lebesgue null set, we may work with $u$ at which
all these inequalities hold, all derivatives involved are defined, and $g_{j,n}(u)\uparrow s_j(u)$ for every
$j$.
Fix such a $u$ and $m\ge0$. For each $n\ge1$ and each $j=0,\ldots,m$, choose $x_j$ from
$\{x_{j,1},\ldots,x_{j,n}\}$ so that $-(a_jx_j)'(u)=g_{j,n}(u)$. Such a choice exists because the maximum of
finitely many values is attained. The derivative inequality for these choices gives $\mu'(u)\ge\sum_{j=0}^m
g_{j,n}(u)$. With $m$ fixed, letting $n\to\infty$ yields $\mu'(u)\ge\sum_{j=0}^m s_j(u)$. This inequality holds
for every $m$ at the same $u$. Since $s_j\ge0$, letting $m\to\infty$ proves the slope inequality in
\eqref{eq_exact_slope}.

Conversely, suppose the slope inequality holds. Fix $m\ge0$ and $x_0,\ldots,x_m\in\KK$, and let
$q_m=\mu+\sum_{j=0}^m a_jx_j$, using the locally absolutely continuous versions in Assumption~\ref{ass_slope}.
The definition of $s_j$ gives $-(a_jx_j)'\le s_j$ \AlmostEverywhere{} for each selected $x_j$. Hence
$q_m'=\mu'+\sum_{j=0}^m(a_jx_j)'\ge\mu'-\sum_{j=0}^m s_j\ge0$ \AlmostEverywhereEnd{} For any $0<u<v<1$, local
absolute continuity gives $q_m(v)-q_m(u)=\int_u^v q_m'(w)\Differential{w}\ge0$. Thus this version of $q_m$ is
nondecreasing on $(0,1)$. Since $q_m\in L^2(0,1)$, it belongs to $\QQ$. This holds for every finite $m$ and
every choice of $x_0,\ldots,x_m$, so condition~\ref{item_finite_support} of Theorem~\ref{thm_support} gives
$X_t\in\QQ$ \AlmostSureEnd{}
\end{proof}

\clearpage
\setcounter{section}{0}\setcounter{equation}{0}\setcounter{example}{0}\setcounter{table}{0}
\renewcommand{\thesection}{S\arabic{section}}\renewcommand{\theequation}{S\arabic{equation}}
\renewcommand{\theproposition}{S\arabic{proposition}}\renewcommand{\theexample}{S\arabic{example}}
\renewcommand{\theremark}{S\arabic{remark}}\renewcommand{\thetable}{S\arabic{table}}
\begin{center}
{\Large\bfseries Supplementary Material\par}
\end{center}
\vspace{1ex}
Results, equations, tables and sections with an S prefix belong to this Supplementary Material; all others refer to the main paper. We use the notation and assumptions of the main paper.

In addition, write $\mathbb{L}^2_{\HH}$ for the space of $\HH$-valued random elements $Y$ on $(\Omega,\mathcal F,\ProbMeasure)$ with $\Expect{\norm{Y}^2}<\infty$, and $\mathbb{L}^2_{\mathbb{R}}$ for the space of real random variables with finite second moment. Under Assumption~\ref{ass_linear}, the series in \eqref{eq_linear} converges in $\mathbb{L}^2_{\HH}$ and \AlmostSure{}: its partial sums form a Hilbert-valued martingale with uniformly bounded second moments, and the tail variances tend to zero.

\section{Geometry and affine restrictions}\label{sec_geometry_s}
This section complements Section~\ref{sec_geometry} of the main paper. Table~\ref{tab_geometry} summarizes the restrictions from Theorem~\ref{thm_collapse} and
Corollary~\ref{cor_limits}. The examples fall into three groups. The first group, the quantile cone $\QQ$ in \eqref{eqquantilecone}, is treated in Section~\ref{sec_geometry} of the main paper and occupies rows 1 and 2 of Table~\ref{tab_geometry}. Section~\ref{sec_geometry} also notes that the cones in the second group are pointed; the details are given below. The third group, convex functions, and the affine restrictions in Proposition~\ref{prop_affine} are treated only here.
\par\medskip
\noindent\textbf{(ii) Nonnegative functions and positive semidefinite operators}\par\nobreak\smallskip
\noindent The nonnegative cone in $L^2(0,1)$ is self-dual and pointed: $f\ge0$ and $-f\ge0$ imply $f=0$
\AlmostEverywhereEnd{} Probability densities additionally satisfy $\int_0^1f(u)\Differential{u}=1$, an affine
restriction treated in Proposition~\ref{prop_affine} below. The consequences of density normalization for
accumulated increments are discussed in Remark~\ref{rem_density} and by \citet{Beare2017}. For covariance
operators, take the Hilbert space of self-adjoint Hilbert--Schmidt operators on a real separable Hilbert space
$\mathcal G$. Write $B\succeq0$ when $B$ is positive semidefinite, and $B\preceq D$ when $D-B\succeq0$. Its
positive semidefinite cone is closed and pointed: $B\succeq0$ and $-B\succeq0$ imply $B=0$. Covariance
operators of square-integrable $\mathcal G$-valued random elements are trace class and therefore belong to this
Hilbert space.
\par\medskip
\noindent\textbf{(iii) Convex functions}\par\nobreak\smallskip
\noindent Let $\CC_{\mathrm{cvx}}$ be the closed cone of functions in $L^2(0,1)$ with a convex version, and
write $\iota(u)=u$. A function and its negative are both convex precisely when the function is affine. Hence
$\lin(\CC_{\mathrm{cvx}})=\spn\{\one,\iota\}$ and $\ncap(\CC_{\mathrm{cvx}})=2$. Under
Theorem~\ref{thm_collapse}, randomness can enter the intercept and slope while curvature remains fixed.
Requiring the functions to be nondecreasing reduces lineality to $\spn\{\one\}$. Adding nonnegativity gives a
pointed cone containing the Lorenz curves. The resulting degeneracy under Theorem~\ref{thm_collapse} does not
use endpoint evaluation as a continuous linear functional on $L^2$. The rank-two increment covariance bound is
attained: let $\{A_t\}$ and $\{B_t\}$ be independent Gaussian random walks started at zero, each with
independent standard normal increments. Then $Q_t(u)=u^2+A_t+B_tu$, $t\ge0$, remains convex. Its centred
independent stationary increments have long-run covariance $\Lambda=\one\otimes\one+\iota\otimes\iota$, which
has rank two because $\one$ and $\iota$ are linearly independent.
\begin{table}[t]
\setlength{\tabcolsep}{3pt}
\renewcommand{\arraystretch}{1.12}
\tbl{Geometry governing the restrictions. The stationary linear process column applies under
Theorem~\ref{thm_collapse}. The increment covariance column applies only when the limit in
Corollary~\ref{cor_limits} exists. Here $c_t$ and $d_t$ are centred scalar linear processes, $\one(u)=1$,
$\iota(u)=u$, and $\Lambda$ is that long-run covariance operator.}{%
\begin{tabularx}{\linewidth}{@{}r@{\enspace}>{\raggedright\arraybackslash}X l l l@{}}
\toprule
& Admissible representation & Lineality & \shortstack[l]{Stationary linear\\process} & \shortstack[l]{Increment\\covariance} \\
\midrule
1 & Real-line quantiles & $\spn\{\one\}$ & $X_t=\mu+c_t\one$ & $\operatorname{rank}\Lambda\le1$ \\
2 & Nonnegative/\allowbreak centred quantiles & $\{0\}$ & $X_t=\mu$ & $\Lambda=0$ \\
\addlinespace[0.35em]
3 & Nonnegative functions (densities) & $\{0\}$ & $X_t=\mu$ & $\Lambda=0$ \\
4 & Positive semidefinite operators & $\{0\}$ & $X_t=\mu$ & $\Lambda=0$ \\
\addlinespace[0.35em]
5 & Convex functions & $\spn\{\one,\iota\}$ & $X_t=\mu+c_t\one+d_t\iota$ & $\operatorname{rank}\Lambda\le2$ \\
6 & Convex nondecreasing functions & $\spn\{\one\}$ & $X_t=\mu+c_t\one$ & $\operatorname{rank}\Lambda\le1$ \\
7 & Nonnegative convex nondecreasing functions (including Lorenz curves) & $\{0\}$ & $X_t=\mu$ & $\Lambda=0$ \\
\bottomrule
\end{tabularx}}
\ReferenceTitle{Geometry governing the restrictions}\label{tab_geometry}
\end{table}
\par\medskip
\noindent\textbf{Affine restrictions}\par\nobreak\smallskip
\noindent Density observations and quantiles with a fixed mean impose affine restrictions in addition to cone
membership. Unit mass fixes the density integral, while a prescribed distributional mean fixes the quantile
integral. The next proposition shows how these restrictions enter admissibility and can remove random
directions left available by the cone.
\begin{proposition}[Affine restrictions within a cone]\label{prop_affine}
Let $p\ge1$ be an integer, let $B:\HH\to\RR^p$ be bounded and linear, let $b\in\RR^p$, and let
$D=\{x\in\CC:Bx=b\}$. Under Assumption~\ref{ass_linear}, $X_t\in D$ \AlmostSure{} if and only if the support
condition in Theorem~\ref{thm_support}\ref{item_support_budget} holds and
\begin{equation}\label{eq_affine}
B\mu=b,\qquad Ba_j\Sigma^{1/2}=0\quad\text{for every }j\ge0.
\end{equation}
Under unbounded lower tails (Definition~\ref{def_unbounded}), $X_t-\mu\in\lin(\CC)\cap\ker B$ \AlmostSureEnd{}
Under the normalized convergence assumption of Theorem~\ref{thm_accumulation}, \eqref{eq_affine} gives
$BF_m\Sigma^{1/2}=0$ for every $m$ and hence $B\widetilde F_\infty\Sigma^{1/2}=0$. Together with
\eqref{eq_accumulation_range}, this yields $\ran(\widetilde F_\infty\Sigma^{1/2})\subseteq\lin(\CC)\cap\ker B$.
\end{proposition}
On $(0,1)$ these restrictions use $Bf=\int_0^1f(u)\Differential{u}$, with $b=1$ for densities and $b$ equal to
the prescribed mean for quantiles. For centred quantiles, $\lin(\QQ)\cap\ker B=\{0\}$, so an additive model with unbounded lower tails is deterministic. For intersections, $\lin(\CC\cap V)=\lin(\CC)\cap V$ for any linear subspace $V$; the
lineality dimension can decrease and never increases. Affine endpoint constraints require a space in which the
relevant traces are defined and continuous, or a direct argument on the chosen versions.

\section{Further results on innovation support and persistence}\label{sec_support_notes}
This section complements Sections~\ref{sec_support} and~\ref{sec_collapse} of the main paper. Section~\ref{sec_support_criterion_s} gives further consequences of Theorem~\ref{thm_support} and discusses the role of independent innovations, Section~\ref{sec_unbounded_s} discusses Definition~\ref{def_unbounded} and Theorem~\ref{thm_collapse}, and Section~\ref{sec_scalar_s} complements Theorem~\ref{thm_accumulation} for scalar filters.

\subsection{Consequences of the support criterion}\label{sec_support_criterion_s}
Since $0\in\KK$, taking $m=0$ and $x_0=0$ in condition~\ref{item_finite_support} of Theorem~\ref{thm_support} gives $\mu\in\CC$. The
characterization also gives a boundary restriction. If $h\in\CC^*$ and $\ip{h}{\mu}=0$,
condition~\ref{item_support_budget} forces every $\kappa(a_j^*h)$ to be zero, since each is nonpositive. Thus
the scalar random variable $\ip{h}{a_j\varepsilon_0}$ is nonnegative \AlmostSure{} and has expectation zero, so
it equals zero \AlmostSure{} for every $j$. Consequently, $\mu\in\lin(\CC)$ forces $X_t-\mu\in\lin(\CC)$ under
any innovation law satisfying Assumption~\ref{ass_linear}.

\begin{remark}[Role of innovation independence]\label{rem_independence}
Independence makes the essential lower bounds of the scalar terms $\ip{h}{a_j\varepsilon_{t-j}}$ additive,
producing the sum $\sum_j\kappa(a_j^*h)$ in condition~\ref{item_support_budget} and making that condition
necessary. If an \AlmostSure{} convergent filtered series is built from dependent innovations whose marginal
supports lie in $\KK$, condition~\ref{item_finite_support} of Theorem~\ref{thm_support} still ensures
admissibility. The converse can fail. For a simple counterexample, take $\HH=\RR$ and $\CC=[0,\infty)$, and let
$\{\eta_t\}$ be independent Rademacher variables. Set $\varepsilon_t=\eta_t-\eta_{t-1}$, $\mu=2$,
$a_0=a_1=\Id$, and $a_j=0$ for $j\ge2$. The sequence $\{\varepsilon_t\}$ is centred and 1-dependent, and
\[
X_t=2+\varepsilon_t+\varepsilon_{t-1}=2+\eta_t-\eta_{t-2}\in\{0,2,4\}
\quad\text{\AlmostSure{}}
\]
Thus the process is admissible and nondegenerate. However, $\KK=[-2,2]$ and $\kappa(1)=-2$, so
condition~\ref{item_support_budget} is violated (at $h=1$, $\ip{h}{\mu}+\sum_{j\ge0}\kappa(a_j^*h)=-2$).
\end{remark}

\subsection{Unbounded lower tails}\label{sec_unbounded_s}
If $\Sigma^{1/2}g=0$, the projection $\ip{g}{\varepsilon_0}$ has variance zero and, by centring, mean zero.
Hence $\ip{g}{\varepsilon_0}=0$ \AlmostSure{} and its essential lower bound $\kappa(g)$ is zero. The definition
thus allows covariance degeneracy. Because it applies to both $g$ and $-g$, every $\ip{g}{\varepsilon_0}$ with
positive variance is also unbounded above. An unbounded innovation norm alone does not imply the condition, nor does the presence of a Gaussian component
that is degenerate in a direction with positive total innovation variance.

\begin{example}[Unbounded support can still be admissible]\label{ex_one_sided}
Let $\HH=\RR$, $\CC=[0,\infty)$, and $\varepsilon_t=E_t-1$ with $E_t$ independent standard exponential
variables. For $a_j\ge0$ with $\sum_j a_j<\infty$, the model is admissible if and only if $\mu\ge\sum_j a_j$.
The innovations are unbounded above but bounded below. This is a nondegenerate process in a pointed cone, so
``unbounded innovations'' alone cannot replace Definition~\ref{def_unbounded}.
\end{example}

\subsection{Scalar filters and persistence}\label{sec_scalar_s}
\begin{remark}[Scalar multiples of a fixed operator]\label{rem_scalar}
Scalar coefficient sequences and scalar multiples of a fixed operator occur in functional local Whittle methods
\citep{LiRobinsonShang2021,SeoShang2024}. If $a_j=b_jB$ for some $B\in\mathcal B(\HH)$ under
Assumption~\ref{ass_linear} and $\sum_j|b_j|=\infty$, admissibility implies
$\ran(B\Sigma^{1/2})\subseteq\lin(\CC)$, even when the signed cumulative coefficients remain bounded. Indeed,
if $Y=\ip{h}{B\varepsilon_0}$ has positive variance for some $h\in\CC^*$, centring implies that $Y$ takes both
signs with positive probability. Hence some $c_h>0$ satisfies $\kappa(b_jB^*h)\le-c_h|b_j|$ for every $j$. When
$\sum_j|b_j|=\infty$, variation in $B\varepsilon_t$ outside lineality allows accumulated negative contributions
that no fixed mean can offset. Admissibility therefore requires $B\varepsilon_t$ to lie in lineality
\AlmostSureEnd{} Since every coefficient is a scalar multiple of $B$, the entire centred process lies there as
well and is zero when the cone is pointed. This applies to the usual scalar fractional filter, whose weights
have order $j^{d-1}$ for $0<d<1/2$. Further details are given in Section~\ref{sec_scalar_details}.
\end{remark}

\section{Further results on processes with stationary increments}\label{sec_increments_notes}

This section complements Section~\ref{sec_increments} of the main paper. Section~\ref{sec_increments_conditions} discusses the conditions of Theorem~\ref{thm_increments}, Section~\ref{sec_increments_corollary} explains how Corollary~\ref{cor_limits} follows from it and when its part~(ii) applies, and Section~\ref{sec_density_normalization} treats random densities whose mean may change over time.

\subsection{Conditions of Theorem~\ref{thm_increments}}\label{sec_increments_conditions}
The centring condition $\Expect{u_1}=0$ in Theorem~\ref{thm_increments} is essential. With $\HH=\RR$ and $\CC=[0,\infty)$, the deterministic increments $u_t=1$ give admissible levels $Q_t=Q_0+t$, although $\lin(\CC)=\{0\}$.

The theorem does not rule out persistence in the levels. If $\{Q_t\}_{t\ge0}$ is any strictly stationary $\CC$-valued process with $\Expect{\norm{Q_0}^2}<\infty$, the increments $u_t=Q_t-Q_{t-1}$ are centred and strictly stationary, and $S_T=Q_T-Q_0$ satisfies $\Expect{\norm{S_T}^2}\le4\Expect{\norm{Q_0}^2}$ however slowly the serial covariance of $Q_t$ decays. Theorem~\ref{thm_increments} restricts the accumulation of increments, not the dependence of the levels.

\subsection{From Theorem~\ref{thm_increments} to Corollary~\ref{cor_limits}}\label{sec_increments_corollary}
For each fixed $h\in\CC^*$, \eqref{eq_uniform_projection} in Theorem~\ref{thm_increments} bounds $\Expect{\ip{h}{S_T}^2}$ uniformly in $T$. Hence $\ip{h}{S_T}/r_T\to0$ in $\mathbb{L}^2_{\mathbb{R}}$ for any $r_T\to\infty$, and $T^{-1}\Expect{\ip{h}{S_T}^2}\to0$. In part~(i) of Corollary~\ref{cor_limits}, the first limit gives $\ip{h}{G}=0$ \AlmostSure{}; in part~(ii), the second gives $\ip{h}{\Lambda h}=0$ and hence $\Lambda h=0$, since $\Lambda$ is nonnegative. Both hold for every $h\in\CC^*$, and \eqref{eq_lineality} then confines $G$ and $\ran\Lambda$ to $\lin(\CC)$. The rank bound in part~(ii) follows because $\ncap(\CC)=\dim\lin(\CC)$ is the noise capacity defined in Section~\ref{sec_support}. The full proof is in Appendix~\ref{app_proofs} of the main paper.

Part~(ii) of Corollary~\ref{cor_limits} assumes only that some bounded operator $\Lambda$ satisfies \eqref{eq_lrv}; it does not assert that these limits exist. They exist, for example, if the lag covariance operators of $\{u_t\}$ are absolutely summable in $\norm{\cdot}_{\mathrm{op}}$; their sum is then a bounded operator $\Lambda$ satisfying \eqref{eq_lrv}.

The argument uses only the scalar bounds for fixed $h\in\CC^*$. It does not require a bound on $\Expect{\norm{PS_T}^2}$, where $P$ is the orthogonal projection onto $\lin(\CC)^\perp$. If $\lin(\CC)^\perp$ is finite dimensional, it is spanned by finitely many elements of $\CC^*$, and \eqref{eq_uniform_projection} gives $\sup_{T\ge1}\Expect{\norm{PS_T}^2}<\infty$. In infinite dimensions, the bounds for each fixed $\ip{h}{S_T}$ do not by themselves imply this norm bound.

\subsection{Random densities with time-varying mean}\label{sec_density_normalization}

\begin{remark}[Density normalization]\label{rem_density}
Let $\{f_t\}_{t\ge0}$ be random probability densities in $L^2(0,1)$ with $\Expect{\norm{f_t}^2}<\infty$, and
put $w_t=f_t-\Expect{f_t}$. The mean $\Expect{f_t}$ may vary with $t$, so the increments $f_t-f_{t-1}$ need not
be centred, and Theorem~\ref{thm_increments} does not in general apply to $Q_t=f_t$. Normalization gives an
\AlmostSureAdjective{} bound instead. For a bounded measurable function $h$ on $(0,1)$, recall that
$\ip{h}{f_t}=\int_0^1h(u)f_t(u)\Differential{u}$. Since $f_t\ge0$ and $\int_0^1f_t(u)\Differential{u}=1$, both
$\ip{h}{f_t}$ and $\ip{h}{\Expect{f_t}}$ lie in $[-\norm{h}_\infty,\norm{h}_\infty]$, so
\[
|\ip{h}{w_T-w_0}|\le4\norm{h}_\infty\quad\text{\AlmostSure{}, for every }T\ge1.
\]
For each bounded measurable $h$, the displayed bound gives $\ip{h}{w_T-w_0}/r_T\to0$ in
$\mathbb{L}^2_{\mathbb{R}}$ whenever $r_T\to\infty$. If $(w_T-w_0)/r_T$ converges weakly in $L^2(0,1)$ to $G$,
it follows that $\ip{h}{G}=0$ \AlmostSure{} for each such $h$. Bounded measurable functions are dense in
$L^2(0,1)$, so $G=0$ \AlmostSureEnd{} If $\Lambda$ satisfies \eqref{eq_lrv} with $S_T=w_T-w_0$, the same bound
gives $\ip{h}{\Lambda h}=0$ for every bounded measurable $h$, and continuity gives $\Lambda=0$. Neither
conclusion requires stationary increments, a linear filter, or a Gaussian limit. \citet{Beare2017} uses density
normalization to show that the linear-increment specification of \citet{ChangKimPark2016} cannot accommodate a
stochastic trend. The bound above yields the same restriction for any square-integrable random densities.
\end{remark}

\section{Further results on quantile processes}\label{sec_quantile_extra}
This section complements Section~\ref{sec_quantile} of the main paper. Section~\ref{sec_slope_interp} interprets the slope condition in Theorem~\ref{thm_slope}, Section~\ref{sec_ar_example} revisits Example~\ref{ex_ar}, and Section~\ref{sec_many_directions} extends that example to many shape directions.

\subsection{Interpretation of the slope condition}\label{sec_slope_interp}
If the distribution whose quantile function is $\mu$ has a positive density $f_\mu$, then
$\mu'(u)=1/f_\mu(\mu(u))$ wherever the inverse-function relation applies. Here $f_\mu$ is the density
corresponding to $\mu$, rather than an average of the densities represented by $X_t$.
Equation~\eqref{eq_exact_slope} permits larger downward derivative perturbations where this reciprocal density
is larger. Such regions need not be confined to the tails.

A flat quantile interval corresponds to a point mass in the represented distribution. If $\mu(u)=c$ on an
interval $I\subset(0,1)$, \eqref{eq_exact_slope} requires $s_j=0$ \AlmostEverywhere{} on $I$ for every $j$. An
admissible $X_t$ must also be constant across $I$, but this value may change with $t$. For example, if
$\{L_t\}$ is a centred square-integrable scalar linear process, $X_t=\mu+L_t\one$ remains a quantile and
equals $c+L_t$ on $I$. Thus the point mass can move in location even though the quantile remains flat on $I$.

The slope criterion in \eqref{eq_exact_slope} relies on local absolute continuity: without it, a nonnegative
derivative \AlmostEverywhere{} need not imply monotonicity. The general support characterization in
Theorem~\ref{thm_support} does not require this assumption.

\subsection{The autoregressive quantile process with one shape direction}\label{sec_ar_example}
Example~\ref{ex_ar} of the main paper takes the uniform quantile function $\mu(u)=u$ and adds one random shape direction $v(u)=u(1-u)$:
\[
X_t(u)=u+Z_tu(1-u),\qquad Z_t=\lambda\sum_{j\ge0}\rho^j\eta_{t-j},
\]
where the $\eta_t$ are independent and identically distributed, centred, with positive variance and $|\eta_t|\le1$, and $\lambda>0$, $|\rho|<1$. Since $v(0)=v(1)=0$, every realization runs from $0$ to $1$, so all represented distributions lie on $[0,1]$. The derivative $X_t'(u)=1+Z_t(1-2u)$ is linear in $u$, with values $1+Z_t$ and $1-Z_t$ at the end points. Hence a realization is nondecreasing exactly when $|Z_t|\le1$. A positive $Z_t$ makes $X_t$ concave, so the represented density $1/X_t'\{X_t^{-1}(x)\}$ increases in $x$ and mass moves towards $1$; a negative $Z_t$ moves mass towards $0$. When $|Z_t|=1$, the slope vanishes at one end point and the density is unbounded there.

Admissibility therefore depends on how large $|Z_t|$ can be. Since $|\eta_t|\le1$, we have $|Z_t|\le\lambda\sum_{j\ge0}|\rho|^j=\lambda/(1-|\rho|)$, so $\lambda\le1-|\rho|$ suffices whatever the law of $\eta_t$. Under the end-point condition $\essinf\eta_0=-1$ and $\esssup\eta_0=1$, this bound is sharp. For each $m$, the event that $\eta_{t-j}$ is close to the sign of $\rho^j$ for $j=0,\ldots,m$ has positive probability, and the remaining terms contribute at most $\lambda|\rho|^{m+1}/(1-|\rho|)$ in absolute value. Hence $\esssup Z_t=\lambda/(1-|\rho|)$, and reversing the signs gives $\essinf Z_t=-\lambda/(1-|\rho|)$. The slope criterion gives the same boundary: here $s_j(u)=\lambda|\rho|^j|1-2u|$, so \eqref{eq_exact_slope} reads $\lambda|1-2u|\le1-|\rho|$ for almost every $u$, and the constraint binds as $u$ approaches $0$ or $1$.

Figure~\ref{fig_boundary} of the main paper plots the two support extremes $q_\pm(u)=u\pm\{\lambda/(1-|\rho|)\}u(1-u)$. Every realization of $X_t$ lies between them, and under the end-point condition realizations come arbitrarily close to each of them with positive probability. They are deterministic bounds, not simulated paths. In the left panel, $\rho=0.6$ and $\lambda=0.3$, so $\lambda/(1-|\rho|)=0.75$. Both extremes are nondecreasing, and so is every realization. In the right panel, $\lambda=1.2$ gives $\lambda/(1-|\rho|)=3$. Then $q_-$ decreases on $(0,1/3)$ and $q_+$ decreases on $(2/3,1)$. Realizations close to either extreme are not quantile functions, and they occur with positive probability, so the model is not admissible. Monotonicity is lost near the end points, where the slope of the shape direction, $|v'(u)|=|1-2u|$, is largest relative to the mean slope $\mu'(u)=1$.

The boundary $\lambda\le1-|\rho|$ links shock size and persistence: for fixed $\lambda$, the extremes widen as $|\rho|$ increases, so more persistent dynamics admit only smaller shocks. Since $\VarOf{Z_t}=\lambda^2\VarOf{\eta_0}/(1-\rho^2)\le\lambda^2/(1-\rho^2)$, admissibility under the end-point condition with $0\le\rho<1$ gives $\VarOf{Z_t}\le(1-\rho)/(1+\rho)$, which tends to zero as $\rho\to1$. The support of $Z_t$ need not fill the interval between its extremes. For Rademacher $\eta_t$, $Z_t/\lambda$ has the law of a symmetric Bernoulli convolution with parameter $|\rho|$, and for $|\rho|<1/2$ its support is a Cantor set \citep{PeresSchlagSolomyak2000}. Only the extremes matter for admissibility.
\subsection{Quantile processes with many shape directions}\label{sec_many_directions}
The process in Section~\ref{sec_ar_example} has a single shape direction. Proposition~\ref{prop_infinite} extends the construction to a countable family of shape directions, using
coefficients $a_j=\lambda\rho^j\Id$ and summable derivative bounds. It shows that admissible shape variation
can have innovation and process covariance operators of infinite rank.
\begin{proposition}[Shape dynamics with bounded innovations]\label{prop_infinite}
Let $\mu$ be as in Assumption~\ref{ass_slope}, with $\mu'\ge c>0$ \AlmostEverywhereEnd{}  Let $v_k$ be
absolutely continuous on $[0,1]$, satisfy $v_k(0)=0$, and obey $V=\sum_{k\ge1}\norm{v_k'}_\infty<\infty$. Let
$\eta_{t,k}$ be independent over $t$ and $k$, identically distributed over $t$ for each $k$, centred, bounded
by one in absolute value, and have positive variance. Set
\begin{equation}\label{eq_infinite_example}
\varepsilon_t=\sum_{k\ge1}\eta_{t,k}v_k,\qquad a_j=\lambda\rho^j\Id,
\qquad \lambda>0,\quad |\rho|<1.
\end{equation}
If $\lambda V\le c(1-|\rho|)$, then the model satisfies Assumption~\ref{ass_linear} and is quantile-valued
\AlmostSureEnd{} If infinitely many $v_k$ are linearly independent, the process covariance has infinite rank.
\end{proposition}
For this construction, $s_j\le\lambda|\rho|^jV$, so Theorem~\ref{thm_slope} gives admissibility under the
stated bound. Example~\ref{ex_ar} is recovered with $\mu(u)=u$, $v_1(u)=u(1-u)$, and $v_k=0$ for $k\ge2$, taking $\eta_{t,1}=\eta_t$.
An infinite-dimensional admissible process can be constructed by taking $\mu(u)=u$, coefficients $\eta_{t,k}$ as in the proposition, for example independent Rademacher variables, and $v_k(u)=2^{-k}(2\pi k)^{-1}\sin(2\pi ku)$. Here $V=1$, so $\lambda\le1-|\rho|$ ensures
admissibility. Since the functions $v_k$ are linearly independent, both the innovation covariance and the
process covariance have infinite rank. Alternatively, the standard normal quantile has $\mu'\ge\sqrt{2\pi}$ and
admits $\lambda V\le\sqrt{2\pi}(1-|\rho|)$. Neither construction automatically concentrates shape variation in
a specific quantile region.

\section{Admissible dynamics for general cones}\label{sec_bounded_s}
This section expands Section~\ref{sec_bounded} of the main paper. We now construct admissible additive dynamics with bounded innovations for a general closed convex cone
$\CC\subseteq\HH$. Recall from Section~\ref{sec_support} that $\KK$ is the closed convex hull of the innovation
support and $a_j$ are the coefficient operators. We give a sufficient condition comparing $|\ip{h}{a_jx}|$ with
$\ip{h}{\mu}$ for $x\in\KK$ and $h\in\CC^*$, and then apply it to densities, covariance matrices, and
quantiles.

Many relevant cones have empty interior in their Hilbert spaces. Recall that $L^2(0,1)$ consists of
square-integrable functions, with $\norm{f}_{L^2}^2=\int_0^1|f(u)|^2\Differential{u}$. A large change on a
sufficiently short interval can have a small $L^2$ norm, yet violate nonnegativity or monotonicity.
Consequently, the nonnegative cone and the quantile cone contain no open $L^2$ ball. For a positive definite
matrix in finite dimensions, sufficiently small perturbations preserve positivity.

For $\mu\in\CC$ and $v\in\HH$, define
\begin{equation}\label{eq_relative_gauge}
p_\mu(v)=\sup_{h\in\CC^*}\frac{|\ip{h}{v}|}{\ip{h}{\mu}},
\end{equation}
where $0/0=0$ and a positive numerator divided by zero is $\infty$. For each $h\in\CC^*$, the ratio compares
the change $|\ip{h}{v}|$ with the mean value $\ip{h}{\mu}$. Applying the dual-cone characterization recalled in
Section~\ref{sec_support} \citep[Section~6.3]{BauschkeCombettes2017} gives
\begin{equation}\label{eq_gauge_geometry}
p_\mu(v)\le1\ \Longleftrightarrow\ \mu+v,\mu-v\in\CC,
\qquad \{v:p_\mu(v)=0\}=\lin(\CC).
\end{equation}
The second identity holds because $p_\mu(v)=0$ is equivalent to all the numerators $|\ip{h}{v}|$ in
\eqref{eq_relative_gauge} being zero, and \eqref{eq_lineality} identifies these $v$ with $\lin(\CC)$. On its
finite-valued domain, $p_\mu$ satisfies $p_\mu(\alpha v)=|\alpha|p_\mu(v)$ for $\alpha\in\RR$ and
$p_\mu(v+w)\le p_\mu(v)+p_\mu(w)$. In particular, $\{v\in\HH:p_\mu(v)<\infty\}$ is a linear subspace. These
properties and \eqref{eq_gauge_geometry} are justified in Section~\ref{sec_relative_details}.

\begin{proposition}[Relative perturbation criterion]\label{prop_gauge}
Suppose Assumption~\ref{ass_linear} holds and $\mu\in\CC$. Put $r_j=\sup_{x\in\KK}p_\mu(a_jx)$. If
$\sum_{j\ge0}r_j\le1$, then $X_t\in\CC$ \AlmostSureEnd{} The conclusion also holds for dependent innovations
satisfying $\varepsilon_t\in\KK$ \AlmostSure{} for every $t$, provided $\mu+\sum_{j=0}^m
a_j\varepsilon_{t-j}\to X_t$ \AlmostSure{} in $\HH$ as $m\to\infty$.
\end{proposition}
This sufficient condition can be conservative because the supremum over $h\in\CC^*$ is taken separately at each
lag. The exact criterion is condition~\ref{item_support_budget} of Theorem~\ref{thm_support}. The mean $\mu$
need not be an interior point of $\CC$: the sufficient condition controls the specified perturbations $a_jx$,
without requiring every sufficiently small perturbation of $\mu$ in the Hilbert norm to remain in $\CC$.
However, not every mean permits nondegenerate variation: if $\mu=0$ and $\CC$ is pointed, admissibility forces
$X_t=0$ \AlmostSure{}, even with bounded innovations.

The following examples apply Proposition~\ref{prop_gauge} to density, covariance matrix, and quantile models.
\begin{example}[Densities]\label{ex_density}
Take $\HH=L^2(0,1)$ and let $\CC$ be its nonnegative cone. Let $\mu\in\CC$ be a density. In this case,
\begin{equation}\label{eq_density_gauge}
p_\mu(v)=\esssup_{u\in(0,1)}\frac{|v(u)|}{\mu(u)},
\end{equation}
with the same zero-denominator convention. The derivation of \eqref{eq_density_gauge} is given in
Section~\ref{sec_density_details}. Let the innovations be independent and identically distributed, centred,
square integrable, have integral zero, and satisfy $|\varepsilon_t(u)|\le R\mu(u)$ \AlmostEverywhere{},
\AlmostSure{}, for some finite $R\ge0$. For $a_j=\lambda\rho^j\Id$, with $\lambda>0$ and $|\rho|<1$,
Proposition~\ref{prop_gauge} gives a density-valued process whenever $\lambda R\le1-|\rho|$.  In particular,
take $\mu(u)=1$ and $\varepsilon_t(u)=\eta_t(2u-1)$ for independent Rademacher variables $\eta_t$. With
$Z_t=\lambda\sum_{j\ge0}\rho^j\eta_{t-j}$ as in Example~\ref{ex_ar}, the model is $X_t(u)=1+Z_t(2u-1)$. It is
density-valued \AlmostSure{} if and only if $\lambda\le1-|\rho|$, and is nondegenerate although the nonnegative
cone is pointed. Thus the assumption of unbounded lower tails in Theorem~\ref{thm_collapse} cannot be omitted.
These claims are verified in Section~\ref{sec_density_details}.
\end{example}
\begin{example}[Covariance matrices]\label{ex_operators}
Let $\HH$ be the space of real symmetric $p\times p$ matrices with inner product $\ip{B}{D}=\TraceOf{BD}$, and
let $\CC$ be its positive semidefinite cone. For a positive definite mean $\mu=M$,
\begin{equation}\label{eq_operator_gauge}
p_M(B)=\norm{M^{-1/2}BM^{-1/2}}_{\mathrm{op}},
\end{equation}
where $\norm{\cdot}_{\mathrm{op}}$ is the operator norm for matrices acting on $\RR^p$.
Let $\{D_t\}$ be centred \IndependentIdentical{} symmetric matrices with $\norm{D_t}_{\mathrm{op}}\le
R<\infty$ \AlmostSureEnd{} Set $\varepsilon_t=M^{1/2}D_tM^{1/2}$ and $a_j=\lambda\rho^j\Id$, with $\lambda>0$
and $|\rho|<1$. Proposition~\ref{prop_gauge} ensures $X_t\succeq0$ \AlmostSure{} when $\lambda R\le1-|\rho|$.
Taking $D_t=\eta_tD$ for independent Rademacher variables $\eta_t$ and a fixed nonzero symmetric matrix $D$
gives a nondegenerate example with $R=\norm{D}_{\mathrm{op}}$. The formula and construction are verified in
Section~\ref{sec_operator_details}.
\end{example}
\begin{example}[Quantiles]\label{ex_quantile_general}
Take $\HH=L^2(0,1)$ and $\CC=\QQ$, the quantile cone in \eqref{eqquantilecone}.  If $\mu\in\QQ$ and $v\in L^2(0,1)$ have
locally absolutely continuous versions, then
\begin{equation}\label{eq_quantile_relative}
p_\mu(v)=\esssup_{u\in(0,1)}\frac{|v'(u)|}{\mu'(u)},
\end{equation}
with the same zero-denominator convention. Constants satisfy $p_\mu(v)=0$. The derivation of
\eqref{eq_quantile_relative} is given in Section~\ref{sec_quantile_details}.
\end{example}
Thus Proposition~\ref{prop_gauge} gives sufficient conditions through derivative bounds for quantiles,
pointwise bounds for densities, and matrix inequalities for covariance matrices.

\section{Proofs}\label{sec_proofs}
This section proves the results stated in this Supplementary Material; the results of the main paper are proved in its Appendix~\ref{app_proofs}. We use the innovation support $\KK=\cconv(\supp\varepsilon_0)$, its lower support functional $\kappa(g)=\inf_{x\in\KK}\ip{g}{x}$, and the innovation covariance $\Sigma$ from Section~\ref{sec_support}.

\begin{proof}[Proof of Proposition~\ref{prop_affine}]
Recall $D=\{x\in\CC:Bx=b\}$, where $B:\HH\to\RR^p$ is bounded and linear. Suppose $BX_t=b$ \AlmostSureEnd{}
Taking expectations gives $B\mu=b$, so $BX_t-B\mu=0$ \AlmostSureEnd{} For each finite $m$, independence and
centring give $\Expect{\ip{Ba_j\varepsilon_{t-j}}{Ba_k\varepsilon_{t-k}}}=0$ whenever $j\ne k$. Expanding the
squared norm therefore yields
\[
\Expect{\norm{\sum_{j=0}^m Ba_j\varepsilon_{t-j}}^2}
=\sum_{j=0}^m\Expect{\norm{Ba_j\varepsilon_0}^2}.
\]
Since $B$ is bounded, the finite sums on the left converge in mean square to $BX_t-B\mu$. Letting $m\to\infty$
gives $0=\Expect{\norm{BX_t-B\mu}^2}=\sum_{j\ge0}\Expect{\norm{Ba_j\varepsilon_0}^2}$. Every summand is
nonnegative, so $\Expect{\norm{Ba_j\varepsilon_0}^2}=0$ and hence $Ba_j\varepsilon_0=0$ \AlmostSure{} for every
$j$.
For $q\in\RR^p$, $\Expect{\ip{q}{Ba_j\varepsilon_0}^2}=\norm{\Sigma^{1/2}a_j^*B^*q}^2$, so this is equivalent
to $Ba_j\Sigma^{1/2}=0$ for every $j$. Conversely, these operator identities give $Ba_j\varepsilon_{t-j}=0$
\AlmostSure{} for every $j$, or equivalently $a_j\varepsilon_{t-j}\in\ker B$. Continuity of $B$ and convergence
of the series then give $BX_t=B\mu=b$. Combining this equivalence with Theorem~\ref{thm_support} proves the
characterization of $X_t\in D$.

Under unbounded lower tails, Theorem~\ref{thm_collapse} gives $X_t-\mu\in\lin(\CC)$ \AlmostSureEnd{} Since
$B(X_t-\mu)=0$, we obtain $X_t-\mu\in\lin(\CC)\cap\ker B$. Finally, \eqref{eq_affine} implies
$BF_m\Sigma^{1/2}=0$ for every $m$. If $r_m^{-1}F_m\to\widetilde F_\infty$ in operator norm, continuity gives
$B\widetilde F_\infty\Sigma^{1/2}=0$. Together with \eqref{eq_accumulation_range}, this proves $\ran(\widetilde
F_\infty\Sigma^{1/2})\subseteq\lin(\CC)\cap\ker B$.
\end{proof}

\begin{proof}[Proof of Proposition~\ref{prop_infinite}]
Recall $V=\sum_{k\ge1}\norm{v_k'}_\infty$ and $a_j=\lambda\rho^j\Id$. Since $v_k$ is absolutely continuous and
$v_k(0)=0$, we have $v_k(u)=\int_0^u v_k'(w)\Differential{w}$. Hence $|v_k(u)|\le
u\norm{v_k'}_\infty\le\norm{v_k'}_\infty$ for $u\in[0,1]$, giving $\norm{v_k}_\infty\le\norm{v_k'}_\infty$.
Since $|\eta_{t,k}|\le1$ and $\norm{v_k}_\infty\le\norm{v_k'}_\infty$, the innovation partial sums satisfy
\[
\sup_{u\in[0,1]}\left|\varepsilon_t(u)-\sum_{k=1}^N\eta_{t,k}v_k(u)\right|
\le\sum_{k>N}\norm{v_k'}_\infty\to 0.
\]
Thus their approximation error tends to zero uniformly over $u$, and
$\norm{\varepsilon_t}\le\norm{\varepsilon_t}_\infty\le V$. Also, $\sum_{k\ge1}|\eta_{t,k}v_k'(w)|\le V$ for
\AlmostEvery{} $w$. Define $d_t(w)=\sum_{k\ge1}\eta_{t,k}v_k'(w)$, so $|d_t(w)|\le V$ \AlmostEverywhereEnd{}
Since $v_k(u)=\int_0^u v_k'(w)\Differential{w}$, dominated convergence gives $\varepsilon_t(u)=\int_0^u
d_t(w)\Differential{w}$. Consequently, $\varepsilon_t$ is absolutely continuous, $\varepsilon_t'=d_t$
\AlmostEverywhere{}, and $\norm{\varepsilon_t'}_\infty\le V$ \AlmostSureEnd{} The assumptions on $\eta_{t,k}$
and dominated convergence show that the innovations are independent and identically distributed and centred.
Moreover,
$\sum_{j\ge0}\Expect{\norm{a_j\varepsilon_0}^2}=\lambda^2\Expect{\norm{\varepsilon_0}^2}/(1-\rho^2)<\infty$,
proving Assumption~\ref{ass_linear}.

For each $m\ge0$, put $q_{t,m}=\mu+\lambda\sum_{j=0}^m\rho^j\varepsilon_{t-j}$. Since $\mu$ is locally
absolutely continuous and each $\varepsilon_{t-j}$ is absolutely continuous, $q_{t,m}$ has a locally absolutely
continuous version. Using $\mu'\ge c$ and $|\varepsilon_{t-j}'(u)|\le V$, we obtain
$q_{t,m}'(u)=\mu'(u)+\lambda\sum_{j=0}^m\rho^j\varepsilon_{t-j}'(u)\ge c-\lambda V\sum_{j=0}^m|\rho|^j\ge
c-\lambda V/(1-|\rho|)\ge0$ \AlmostEverywhere{}, \AlmostSureEnd{} This version is therefore nondecreasing, so
$q_{t,m}\in\QQ$ \AlmostSureEnd{} Assumption~\ref{ass_linear} gives $q_{t,m}\to X_t$ \AlmostSure{} in
$L^2(0,1)$. Since $\QQ$ is closed, $X_t\in\QQ$ \AlmostSureEnd{}

Independence and centring of the coefficients give $\ip{g}{\Sigma
g}=\sum_{k\ge1}\VarOf{\eta_{0,k}}\ip{g}{v_k}^2$ for every $g\in L^2(0,1)$. Since every $\VarOf{\eta_{0,k}}$ is
positive, this sum is zero precisely when $\ip{g}{v_k}=0$ for every $k$. Hence
$\ker\Sigma=(\overline{\spn}\{v_k:k\ge1\})^\perp$. If infinitely many $v_k$ are linearly independent, $\Sigma$
therefore has infinite rank. The process covariance is $\lambda^2(1-\rho^2)^{-1}\Sigma$, so it has the same
rank.
\end{proof}
\begin{proof}[Proof of Proposition~\ref{prop_gauge}]
Recall $r_j=\sup_{x\in\KK}p_\mu(a_jx)$. Since $\sum_{j\ge0}r_j\le1$, each $r_j$ is finite. Fix $j\ge0$,
$x\in\KK$, and $h\in\CC^*$. The definitions give $p_\mu(a_jx)\le r_j$. If $\ip{h}{\mu}>0$, the ratio in
\eqref{eq_relative_gauge} is therefore at most $r_j$, so $|\ip{h}{a_jx}|\le r_j\ip{h}{\mu}$. If
$\ip{h}{\mu}=0$, finiteness of $p_\mu(a_jx)$ forces $\ip{h}{a_jx}=0$ by the zero-denominator convention, giving
$|\ip{h}{a_jx}|\le r_j\ip{h}{\mu}=0$.

For every $m\ge0$ and every $x_0,\ldots,x_m\in\KK$, it follows that
\[
\ip{h}{\mu+\sum_{j=0}^m a_jx_j}
\ge\ip{h}{\mu}-\sum_{j=0}^m|\ip{h}{a_jx_j}|
\ge\Bigl(1-\sum_{j=0}^m r_j\Bigr)\ip{h}{\mu}\ge0.
\]
This holds for every $h\in\CC^*$, so the dual-cone characterization in Section~\ref{sec_support} gives
$\mu+\sum_{j=0}^m a_jx_j\in\CC$. This is condition~\ref{item_finite_support} of Theorem~\ref{thm_support},
which gives $X_t\in\CC$ \AlmostSureEnd{}

The conclusion $\mu+\sum_{j=0}^m a_jx_j\in\CC$ for arbitrary $x_0,\ldots,x_m\in\KK$ was derived without using
independence. For dependent innovations satisfying $\varepsilon_s\in\KK$ \AlmostSure{} for every $s$,
substituting $x_j=\varepsilon_{t-j}$ therefore gives $\mu+\sum_{j=0}^m a_j\varepsilon_{t-j}\in\CC$
\AlmostSure{} for every $m$. When independence is dropped, the second-moment summability condition in
Assumption~\ref{ass_linear} alone need not imply convergence of these partial sums. Their
\AlmostSureAdjective{} convergence to $X_t$ in $\HH$ is therefore assumed separately in the dependent case. On
the probability-one event where the partial sums belong to $\CC$ for every $m$ and converge in $\HH$ to $X_t$,
closedness of $\CC$ implies $X_t\in\CC$.
\end{proof}
\section{Auxiliary results and complementary discussion}\label{sec_auxiliary}
This section gives the details for Remark~\ref{rem_scalar} and the derivations for Section~\ref{sec_bounded_s}. The countable representation used in the proof of Theorem~\ref{thm_slope} is Lemma~\ref{lem_essential_family} of the main paper.
\subsection{Scalar multiples of a fixed operator (Remark~\ref{rem_scalar})}\ReferenceTitle{Scalar multiples of
a fixed operator}\label{sec_scalar_details}
For the model $a_j=b_jB$ in Remark~\ref{rem_scalar}, we explain why $\sum_{j\ge0}|b_j|=\infty$ forces
$\ran(B\Sigma^{1/2})\subseteq\lin(\CC)$ under admissibility. The argument does not require the signed sums
$\sum_{j=0}^m b_j$ to diverge.

Fix $h\in\CC^*$ and put $Y=\ip{h}{B\varepsilon_0}$. If $\VarOf{Y}>0$, centring implies that $Y$ takes positive
and negative values with positive probability. Choose $c_h>0$ such that $\ProbOf{Y<-c_h}>0$ and
$\ProbOf{Y>c_h}>0$. On $\{Y<-c_h\}$ when $b_j>0$, and on $\{Y>c_h\}$ when $b_j<0$, we have $b_jY<-c_h|b_j|$.
Hence $\kappa(b_jB^*h)=\essinf(b_jY)\le-c_h|b_j|$ for every $j$, with equality at $b_j=0$. Since
$\sum_{j\ge0}|b_j|=\infty$, it follows that $\sum_{j\ge0}\kappa(b_jB^*h)=-\infty$, contradicting
condition~\ref{item_support_budget} of Theorem~\ref{thm_support}.

Therefore $\VarOf{Y}=\norm{\Sigma^{1/2}B^*h}^2=0$ for every $h\in\CC^*$. Equation~\eqref{eq_lineality},
continuity, and taking adjoints give $\ran(B\Sigma^{1/2})\subseteq\lin(\CC)$. Since $Y$ is centred, we also
have $\ip{h}{B\varepsilon_0}=0$ \AlmostSure{} for each $h\in\CC^*$. Choosing a countable dense subset of
$\CC^*$ and extending by continuity yields $B\varepsilon_0\in(\CC^*)^\perp=\lin(\CC)$ \AlmostSureEnd{} Thus
every $a_j\varepsilon_{t-j}=b_jB\varepsilon_{t-j}$ belongs to $\lin(\CC)$ \AlmostSure{}, as does their
convergent sum. If $\CC$ is pointed, this gives $X_t=\mu$ \AlmostSureEnd{}
\subsection{Formulas and constructions for general cones}\label{sec_relative_details}
We justify the properties of $p_\mu$ and the formulas and constructions in Examples~\ref{ex_density},
\ref{ex_operators}, and~\ref{ex_quantile_general}. Recall \eqref{eq_relative_gauge} and its zero-denominator
convention.

For $r\ge0$, the definition gives $p_\mu(v)\le r$ if and only if $|\ip{h}{v}|\le r\ip{h}{\mu}$ for every
$h\in\CC^*$. When $\ip{h}{\mu}=0$, both conditions require $\ip{h}{v}=0$. By the dual-cone characterization in
Section~\ref{sec_support}, these inequalities are equivalent to $r\mu+v,r\mu-v\in\CC$. Taking $r=1$ proves the
first identity in \eqref{eq_gauge_geometry}; taking $r=0$ gives $p_\mu(v)=0$ if and only if
$v\in(\CC^*)^\perp=\lin(\CC)$, by \eqref{eq_lineality}.

For $v,w$ with finite $p_\mu$ and $\alpha\in\RR$, the definition and the scalar triangle inequality give
$p_\mu(\alpha v)=|\alpha|p_\mu(v)$ and $p_\mu(v+w)\le p_\mu(v)+p_\mu(w)$. Hence $\{v:p_\mu(v)<\infty\}$ is a
linear subspace.

The density and covariance matrix constructions both use $a_j=\lambda\rho^j\Id$, with $\lambda>0$ and
$|\rho|<1$. For centred \IndependentIdentical{} square-integrable innovations, Assumption~\ref{ass_linear}
follows from
$\sum_{j\ge0}\Expect{\norm{a_j\varepsilon_0}^2}=\lambda^2\Expect{\norm{\varepsilon_0}^2}/(1-\rho^2)<\infty$. If
$p_\mu(\varepsilon_0)\le R$ \AlmostSure{} for some finite $R\ge0$, the closed convex set
$\{v:R\mu+v,R\mu-v\in\CC\}$ contains $\KK$. Thus $p_\mu(x)\le R$ for $x\in\KK$, giving $r_j\le\lambda|\rho|^jR$
and $\sum_jr_j\le\lambda R/(1-|\rho|)$. Proposition~\ref{prop_gauge} therefore applies when $\lambda
R\le1-|\rho|$.
\AppendixPart{Densities}\label{sec_density_details}
For the nonnegative cone in $L^2(0,1)$, $r\mu+v,r\mu-v\in\CC$ means that $r\mu(u)+v(u)\ge0$ and
$r\mu(u)-v(u)\ge0$ \AlmostEverywhereEnd{} These inequalities are equivalent to $-r\mu(u)\le v(u)\le r\mu(u)$
\AlmostEverywhereEnd{} Thus the general characterization gives $p_\mu(v)\le r$ if and only if $|v|\le r\mu$
\AlmostEverywhereEnd{} Taking the infimum over $r\ge0$ gives \eqref{eq_density_gauge}. On $\{\mu=0\}$, any
finite bound requires $v=0$ \AlmostEverywhere{}, as encoded by the denominator convention.

The bound $|\varepsilon_t|\le R\mu$ in Example~\ref{ex_density} gives $p_\mu(\varepsilon_t)\le R$
\AlmostSure{}, so Proposition~\ref{prop_gauge} yields $X_t\ge0$ \AlmostSure{} when $\lambda R\le1-|\rho|$. Put
$q_{t,m}=\mu+\sum_{j=0}^m a_j\varepsilon_{t-j}$. Since $\mu$ has integral one and each innovation has integral
zero, $\int_0^1q_{t,m}(u)\Differential{u}=1$ \AlmostSure{} for every $m$. Assumption~\ref{ass_linear} gives
$q_{t,m}\to X_t$ \AlmostSure{} in $L^2(0,1)$. Thus Cauchy--Schwarz gives
$|\int_0^1X_t(u)\Differential{u}-1|=|\int_0^1(X_t-q_{t,m})(u)\Differential{u}|\le\norm{X_t-q_{t,m}}\to0$
\AlmostSureEnd{} Since the left-hand side does not depend on $m$, $\int_0^1X_t(u)\Differential{u}=1$
\AlmostSureEnd{}

In the Rademacher construction of Example~\ref{ex_density}, $1+Z_t(2u-1)$ is nonnegative \AlmostEverywhere{}
exactly when $|Z_t|\le1$. As in Example~\ref{ex_ar}, $\esssup|Z_t|=\lambda/(1-|\rho|)$, giving the exact
condition $\lambda\le1-|\rho|$. Nondegeneracy follows from $\VarOf{Z_t}=\lambda^2/(1-\rho^2)>0$.
\AppendixPart{Covariance matrices}\label{sec_operator_details}
For the positive definite $M$ in Example~\ref{ex_operators}, the general characterization gives $p_M(B)\le r$
if and only if $-rM\preceq B\preceq rM$. Multiplying on the left and right by $M^{-1/2}$ gives $-rI_p\preceq
M^{-1/2}BM^{-1/2}\preceq rI_p$, where $I_p$ is the $p\times p$ identity matrix. Since $M^{-1/2}BM^{-1/2}$ is
symmetric, these inequalities are equivalent to $\norm{M^{-1/2}BM^{-1/2}}_{\mathrm{op}}\le r$, proving
\eqref{eq_operator_gauge}.

For the construction, $\varepsilon_t=M^{1/2}D_tM^{1/2}$ is centred and \IndependentIdentical{}, with
$\norm{\varepsilon_t}\le\sqrt p\,\norm{M}_{\mathrm{op}}R$ and $p_M(\varepsilon_t)=\norm{D_t}_{\mathrm{op}}\le
R$ \AlmostSureEnd{} The coefficient summability verified at the start of Section~\ref{sec_relative_details}
and Proposition~\ref{prop_gauge} give admissibility when $\lambda R\le1-|\rho|$. For $D_t=\eta_tD$ with
$D\ne0$, independence and centring give
$\Expect{\norm{X_t-M}^2}=\lambda^2\norm{M^{1/2}DM^{1/2}}^2/(1-\rho^2)>0$, proving nondegeneracy.
\AppendixPart{Quantiles}\label{sec_quantile_details}
For the locally absolutely continuous $\mu$ and $v$ in Example~\ref{ex_quantile_general}, $p_\mu(v)\le r$ is
equivalent to $r\mu+v,r\mu-v\in\QQ$, or $|v'|\le r\mu'$ \AlmostEverywhereEnd{} Taking the infimum over $r\ge0$
proves \eqref{eq_quantile_relative}. Moreover, $p_\mu(v)=0$ if and only if $v'=0$ \AlmostEverywhere{}, which by
integration on closed subintervals is equivalent to $v$ being constant.

\setlength{\bibsep}{2pt}

\end{document}